\documentclass[a4paper]{amsart}

\usepackage{graphicx}
\newtheorem{theorem}{Theorem}[section]
\newtheorem{proposition}[theorem]{Proposition}
\newtheorem{lemma}[theorem]{Lemma}
\newtheorem{claim}[theorem]{Claim}

\newtheorem{corollary}[theorem]{Corollary}

\theoremstyle{definition}

\theoremstyle{remark}

\numberwithin{equation}{section}

\begin{document}

\title{Triple crossing numbers of graphs}
\author{Hiroyuki Tanaka}
\address{Graduate School of Education,
Hiroshima University, 1-1-1 Kagamiyama, Higashi-hiroshima 739-8524, Japan}
\author{Masakazu Teragaito}
\address{Department of Mathematics and Mathematics Education,
Graduate School of Education,
Hiroshima University, 1-1-1 Kagamiyama, Higashi-hiroshima 739-8524, Japan}
\email{teragai@hiroshima-u.ac.jp}



\begin{abstract}
We introduce the triple crossing number, a variation of the crossing number, of a graph, which
is the minimal number of crossing points in all drawings of the graph 
with only triple crossings.
It is defined to be zero for planar graphs, and to be infinite
for non-planar graphs which  do not admit a drawing with only
triple crossings.
In this paper, we determine the triple crossing numbers for all complete multipartite graphs
which include all complete graphs.
\end{abstract}

\maketitle

\section{Introduction}

Let $G$ be a graph.
A \textit{drawing} of $G$ means a representation
of the graph in the Euclidean plane or the $2$-sphere, where
vertices are points and edges are simple arcs joining their end-vertices.
Since each edge is simple, no edge admits self crossings.
Furthermore, we assume that
the interiors of edges do not contain vertices, and that
two edges do not intersect if they have a common vertex, and 
that two edges without common end-vertex intersect at most once, and if so, then they
intersect transversally.
These requirements are essential in this paper.
A drawing is called a \textit{regular drawing\/} (resp. \textit{semi-regular drawing}\/) if it has only double (resp. triple) crossing points.
From the requirements, we know that a graph has at least $6$ vertices if it admits a semi-regular drawing
with at least one triple crossing point.

The crossing number $\mathrm{cr}(G)$ of $G$ is defined to be
the minimal number of crossing points over all regular drawings of $G$.
In particular, $\mathrm{cr}(G)=0$ if $G$ is planar.
In this paper, we introduce a new variation of the crossing number.
The \textit{triple crossing number\/} $\mathrm{tcr}(G)$ is zero if $G$ is planar,
and $\infty$ if $G$ does not admit a semi-regular drawing.
Otherwise, 
$\mathrm{tcr}(G)$ is defined to be
the minimal number of triple crossing points over all semi-regular drawings of $G$.
In particular, $\mathrm{tcr}(G)=0$ if and only if $G$ is planar.

\begin{figure}[tb]
\includegraphics*[scale=0.7]{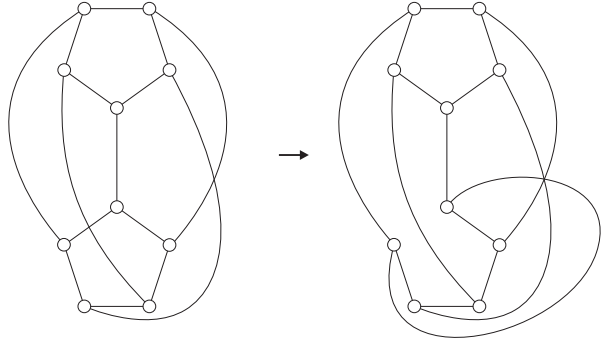}
\caption{The Petersen graph}
\label{fig:petersen}
\end{figure}

The triple crossing number can be regarded as a specialization of the degenerate crossing
number introduced by Pach and T\'{o}th \cite{PT}.
In addition, for example,
the Petersen graph is known to have the crossing number two (and thus non-planar),
and hence has the triple crossing number one from Figure \ref{fig:petersen}.
In general, we have the inequality $\mathrm{cr}(G)\le 3\,\mathrm{tcr}(G)$ for these two notions,
since we obtain a regular drawing from a semi-regular drawing by
perturbing each triple crossing point into three double crossing points.

In this paper, we determine the triple crossing numbers for all complete multipartite graphs.
A complete multipartite graph is a graph whose vertex set
can be partitioned into at least two, mutually disjoint non-empty sets, called
the partite sets,
so that two vertices $u$ and $v$ are adjacent if and only if $u$ and $v$ belong
to different sets of the partition.
If the partite sets are of sizes $n_1,\dots,n_t$\ $(n_i\ge 1)$, then
the graph is denoted by $K_{n_1,\dots,n_t}$.
We always assume that $n_i\ge n_j$ if $i<j$.
In particular, if all $n_i=1$, then the graph $K_{1,\dots,1}$ is the complete graph $K_t$ with $t$ vertices.

Here is how the paper is organized.
After we describe basic lemmas, used in the paper repeatedly, in Section \ref{sec:pre},
we show that the triple crossing number of a complete $t$-partite graph
is $\infty$ if $t\ge 5$
in Section \ref{sec:5partite}.
In the successive sections, we work on the cases when $t\le 4$.
Here we should mention that
the hardest part is the case where $t=2$, in particular,
long, but elementary, geometric arguments are needed
to show that $K_{5,4}$, $K_{4,4}$, $K_{5,3}$ and $K_{n,3}$ with $n\ge 7$
do not admit a semi-regular drawing.
This is treated in Sections \ref{sec:k54}, \ref{sec:k44} and \ref{sec:kn3}.
After concluding the case where $t=2$ in Section \ref{sec:bipartite},
the cases where $t=4$ and $t=3$ are established in Sections \ref{sec:4partite} and
\ref{sec:tripartite}, respectively.
Section \ref{sec:comment} contains some remarks on our requirements for drawings and
a generalization of triple crossing number.

We would like to thank the referee for careful reading and suggestions making the paper
more readable.

\section{Basic lemmas}\label{sec:pre}

Basic terms of graph theory can be found in textbooks such as
\cite{BM, We}.

\begin{lemma}\label{lem:nonplanar}
The complete bipartite graph $K_{3,3}$ and the complete graph $K_5$ with five vertices
are non-planar.
Also, a graph is non-planar if it contains $K_{3,3}$ or $K_5$ as a subgraph.
\end{lemma}

\begin{lemma}\label{lem:plane}
Let $G$ be a plane graph with $p\ (\ge 3)$ vertices and $q$ edges.
Then $G$ has a vertex of degree less than $6$
and the inequality $q\le 3p-6$ holds.
Furthermore, $q=3p-6$ if and only if each region of $G$ is $3$-sided.
\end{lemma}

For the proofs of these lemmas, see \cite{BM}.

\begin{lemma}\label{lem:alg-condition}
Let $G$ be a graph with $p\ (\ge 3)$ vertices and $q$ edges.
If $G$ admits a semi-regular drawing, then $q\le 3p-6$.
Thus, if $q>3p-6$, then $\mathrm{tcr}(G)=\infty$.
\end{lemma}

\begin{proof}
Let $D$ be a semi-regular drawing of $G$, and let $k$ be the number of
triple crossing points in $D$.
If a new vertex is added to each triple crossing point, then we obtain
a (simple) plane graph $G'$.
Since $G'$ has $p+k$ vertices and $q+3k$ edges, we have $q+3k\le 3(p+k)-6$ by Lemma \ref{lem:plane},
from which we have the conclusion.
\end{proof}

\begin{lemma}\label{lem:faces}
Let $G$ be a connected plane graph with $p\ (\ge 3)$ vertices, $q$ edges and $r$ faces. 
Let $d=3p-q-6$. 
\begin{enumerate}
\item If $d=1$, then one face is $4$-sided, and the others are $3$-sided.
\item If $d=2$, then either
 \begin{enumerate}
 \item one face is $5$-sided, and the others are $3$-sided; or
 \item two faces are $4$-sided, and the others are $3$-sided.
 \end{enumerate}
\item If $d=3$, then either
 \begin{enumerate}
 \item one face is $6$-sided, and the others are $3$-sided; or
 \item one face is $5$-sided, another face is $4$-sided, and the others are $3$-sided; or
 \item three faces are $4$-sided, and the others are $3$-sided.
 \end{enumerate}
\end{enumerate}
\end{lemma}

We remark that $d\ge 0$, and $d=0$ if and only if all faces are $3$-sided by Lemma \ref{lem:plane}.

\begin{proof}
By Euler's formula (\cite{BM}), $p-q+r=2$.
Let $r_i$ denote the number of $i$-sided faces of $G$. Then
\begin{equation}\label{eq:eq}
3r_3+4r_4+5r_5+6r_6+\sum_{i\ge 7}ir_i=2q.
\end{equation}
Thus
$7r-4r_3-3r_4-2r_5-r_6\le 2q$.
Since $q=p+r-2$ and $d=2p-r-4$, we have
\begin{equation}\label{eq:diff}
4r-d\le 4r_3+3r_4+2r_5+r_6\le 4(r_3+r_4+r_5+r_6)\le 4r.
\end{equation}
In particular, the difference between the second and third terms, which is
$r_4+2r_5+3r_6$, is at most $d$.
We remark that $d\equiv r\pmod{2}$.

When $d\in \{1,2,3\}$, $4r$ is the only multiple of four within the interval $[4r-d, 4r]$.
Since $4(r_3+r_4+r_5+r_6)$ is a multiple of four, we see
$4(r_3+r_4+r_5+r_6)=4r$, giving $r_3+r_4+r_5+r_6=r$.
Furthermore, if $r_5=r_6=0$,
then (\ref{eq:eq}) reduces to $3r_3+4r_4=2q$.
Combining this with $r_3+r_4=r$ gives $r_4=2q-3r=d$.

(1) Since $r_4+2r_5+3r_6\le d=1$, we have $r_5=r_6=0$.
Then $r_4=1$, and thus $r_3=r-1$.

(2) Since $r_4+2r_5+3r_6\le 2$, we have $r_6=0$ and $r_5\le 1$.
If $r_5=1$, then $r_4=0$, giving $r_3=r-1$.
This is the conclusion (a).
If $r_5=0$, then $r_4=2$, and thus $r_3=r-2$.
This is the conclusion (b).

(3) Since $r_4+2r_5+3r_6\le 3$, we have $r_6\le 1$.
If $r_6=1$, then $r_4=r_5=0$, and thus $r_3=r-1$.
This is the conclusion (a).

Suppose $r_6=0$.
Since $r_4+2r_5\le 3$, we see $r_5\le 1$.

If $r_5=1$, then $r_4\le 1$.
From (\ref{eq:eq}), $3r_3+4r_4=2q-5$.
Combining this with $r_3+r_4=r-1$ gives $r_4=2q-3r-2$.
Since $r\equiv 1\pmod{2}$, we have $r_4=1$.
Hence $r_4=r_5=1$ and $r_3=r-2$.
This is the conclusion (b).

Finally suppose $r_5=0$.
Then $r_4=3$, and thus $r_3=r-3$.
This is the conclusion (c).
\end{proof}

\section{Complete $t$-partite graphs $(t\ge 5)$}\label{sec:5partite}

\begin{theorem}\label{thm:5multi}
If $t\ge 5$, then no complete $t$-partite graph $G$ admits
a semi-regular drawing.
Thus, $\mathrm{tcr}(G)=\infty$.
\end{theorem}

\begin{proof}
Assume for a contradiction that $G$ admits a semi-regular drawing $D$.
Let $t\ge 7$.
If a new vertex is added to each triple crossing point, then
we have a plane graph $G'$.
However, the original vertices have degree at least $t-1\,(\ge 6)$, and the new vertices
have degree $6$.
This contradicts Lemma \ref{lem:plane}.

Let $G=K_{n_1,n_2,\dots,n_6}$.
Then $G$ has $p=\sum_{i}n_i$ vertices
and $q=\sum_{i<j}n_in_j$ edges.
Then
\begin{eqnarray*}
q-3p+6 &=& (n_1+n_4-3)(n_2+n_3-3)+n_1n_4+n_2n_3\\
      &\quad & +(n_5+n_6)(n_1+n_2+n_3+n_4-3)+n_5n_6-3\\
 &\ge & (2n_4-3)^2+2n_4^2\ge 3.
\end{eqnarray*}
This contradicts Lemma \ref{lem:alg-condition}.

Finally, let $G=K_{n_1,n_2,\dots,n_5}$.
As above,
\begin{eqnarray*}
q-3p+6 &=& (n_1+n_4-3)(n_2+n_3-3)+n_1n_4+n_2n_3\\
  &\quad & +n_5(n_1+n_2+n_3+n_4-3)-3\\
 &\ge&  (2n_4-3)^2+2n_4^2+n_5-3\ge 1.
\end{eqnarray*}
This contradicts Lemma \ref{lem:alg-condition} again.
\end{proof}

\begin{corollary}\label{cor:complete}
Let $K_n$ be the complete graph with $n$ vertices.
Then
\[
\mathrm{tcr}(K_n)=
\begin{cases}
0 & \text{if $n\le 4$},\\
\infty & \text{otherwise}.
\end{cases}
\]
\end{corollary}

\begin{proof}
If $n\le 4$, then $K_n$ is planar, and thus $\mathrm{tcr}(K_n)=0$ by definition.
The rest follows from Theorem \ref{thm:5multi}.
\end{proof}

\section{$K_{5,4}$}\label{sec:k54}

Throughout this section, we will assume that $G=K_{5,4}$
admits a semi-regular drawing.
We will show that this is impossible.

Let $V_1=\{x_1,x_2,x_3,x_4,x_5\}$ and $V_2=\{A, B, C, D\}$ be the partite sets of $G$.
For convenience, we refer to vertices of $V_1$ (resp. $V_2$) as \textit{black\/} (resp. \textit{white\/}) vertices.
We denote the edge $Ax_i$ by $a_i$ for $1\le i\le 5$.
These are called \textit{$A$-lines}.
Similarly, define $b_i,c_i,d_i$, and call them $B$-, $C$-, $D$-lines, respectively.
In particular, each black vertex
is of degree four, and is incident with four distinct classes of lines.

We fix a semi-regular drawing of $G$ hereafter, which is denoted by the same symbol $G$.
Notice that any line in $G$ connects a white vertex and a black vertex,
and there may be triple crossing points on it.
From our requirements for drawings, each triple crossing point of $G$
arises from three distinct classes of lines.
This fact will be referred to as \textit{property} ($\ast$) throughout the paper.
Property ($\ast$) is very useful and powerful.
For example, if an $A$-line and a $B$-line intersect at a triple crossing point,
then we can conclude that the remaining line through the triple crossing point
is either a $C$- or $D$-line.

Let $k$ be the number of triple crossing points.
Add a new vertex to each triple crossing point.
Then we have a plane graph $G'$ with $9+k$ vertices and $20+3k$ edges.
Since $3(9+k)-(20+3k)-6=1$, 
the faces of the plane graph $G'$ are all $3$-sided, except a single $4$-sided face
by Lemma \ref{lem:faces}.
For the semi-regular drawing $G$,
a face means that of $G'$, although it is an abuse of words.
A $3$-sided face is also called a triangle.

\begin{figure}[tb]
\includegraphics*[scale=0.7]{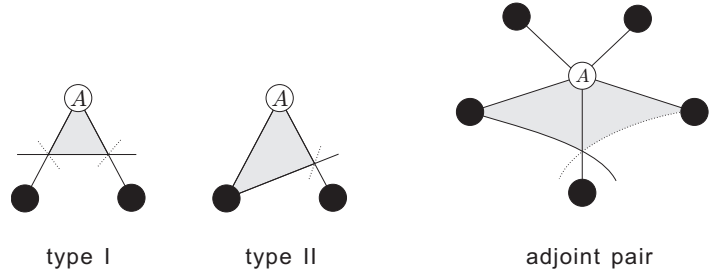}
\caption{Two types of triangles at $A$ and an adjoint pair of type II triangles}
\label{fig:triangle}
\end{figure}

Take a look around vertex $A$.
There are five faces of $G$, since $A$ is not a cut vertex.
We may assume that all five faces around vertex $A$ are triangles without loss of generality,
since the $4$-sided face is incident with
at most two white vertices.
There are two types of triangle around vertex $A$
as shown in Figure \ref{fig:triangle}.
A \textit{type I triangle\/} is incident with
two triple crossing points, and
a \textit{type II triangle\/} is incident with
a black vertex and a triple crossing point.

Notice that type II triangles appear in pairs.
More precisely, this means that 
every type II triangle at $A$ shares
an $A$-line fully with another type II triangle.
Such a pair of type II triangle is referred to as an \textit{adjoint pair\/} of type II triangles.
See Figure \ref{fig:triangle}.
Hence the number of type II triangles at vertex $A$ is either $0$, $2$ or $4$.
We will eliminate these three possibilities.

\begin{lemma}\label{lem:typeII-4}
The number of type II triangles at vertex $A$ is not four.
\end{lemma}

\begin{proof}
Suppose that there are four type II triangles at $A$.
Then 
we can assume that the local configuration at $A$ is
as shown in Figure \ref{fig:4triangle}(1) by renaming $B, C, D$, if necessary.
(Recall that each black vertex is incident with four distinct classes of lines.)
Then, by property ($\ast$),
the horizontal line is a $B$- or $D$-line, and
the right lower line is a $C$- or $D$-line.
See Figure \ref{fig:4triangle}(2), 
where the symbol $B/D$, for example,  indicates the class of the horizontal line.
It does not mean that vertex $B$ or $D$ locates the left side of the horizontal line.
Thus there are four cases as shown in Figure \ref{fig:case1-sub}.
Here, the class of the right upper line is determined by property ($\ast$) and the fact
that each black vertex is incident with four distinct classes of lines.
For Figure \ref{fig:case1-sub}(4), the right upper line may be a $C$-line.
But then, it can be reduced to (3) by renaming $B$ and $D$ and symmetry.

\begin{figure}[tb]
\includegraphics*[scale=0.7]{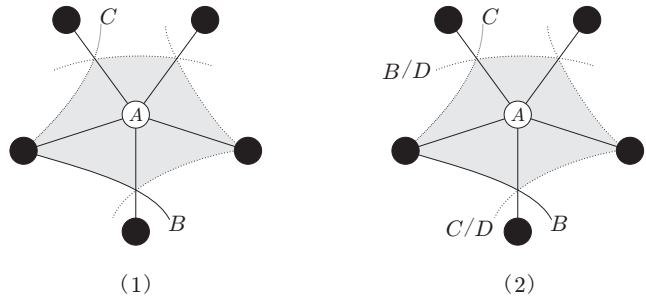}
\caption{Four type II triangles at $A$}
\label{fig:4triangle}
\end{figure}

\begin{figure}[tb]
\includegraphics*[scale=0.7]{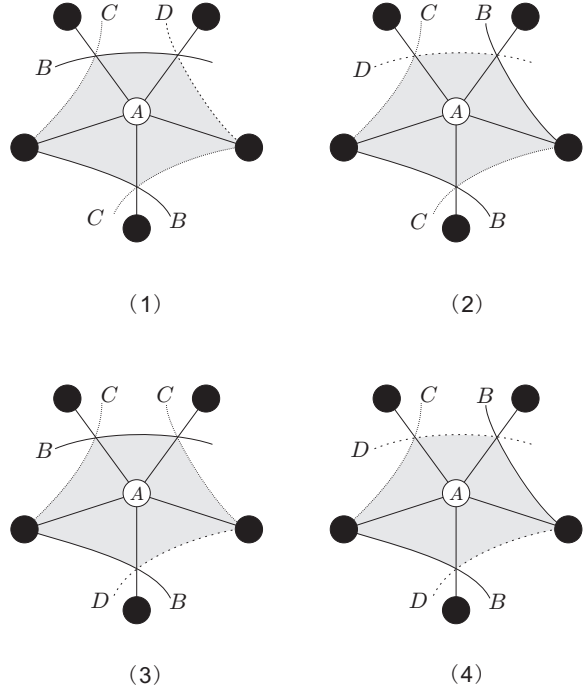} 
\caption{Four cases where there are four type II triangles}
\label{fig:case1-sub}
\end{figure}

\begin{claim}\label{cl:case1}
Figure \ref{fig:case1-sub}(1) is impossible.
\end{claim}

\begin{proof}
Consider the $B$-line $b_2$.
Let $f_1$ and $f_2$ be the faces
incident with $b_2$ and vertex $x_2$.
See Figure \ref{fig:case1-1-2}.

\begin{figure}[tb]
\includegraphics*[scale=0.7]{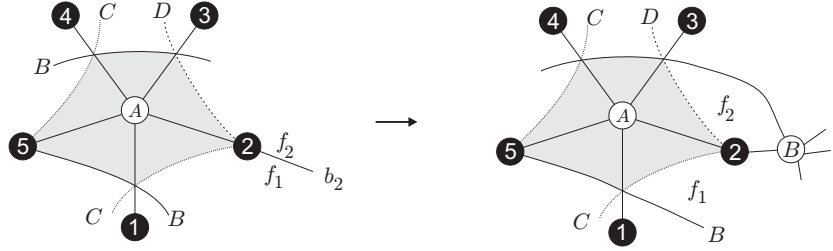}
\caption{The $B$-line $b_2$ and two faces $f_1$, $f_2$}
\label{fig:case1-1-2}
\end{figure}

If $f_1$ is $4$-sided, then $f_2$ is $3$-sided, because there is only one
$4$-sided face.
Since two $B$-lines cannot intersect at a triple crossing point by property ($\ast$),
the horizontal $B$-line and $b_2$ meet at vertex $B$
as shown in the second of Figure \ref{fig:case1-1-2}
in order to make $f_2$ $3$-sided.
However, then $f_1$ cannot be $4$-sided.
Thus we can conclude that $f_1$ is $3$-sided.
Then vertex $B$ is located around $f_1$ as shown in the first of Figure \ref{fig:case1-1}
by the same reason as above.
\begin{figure}[tb]
\includegraphics*[scale=0.7]{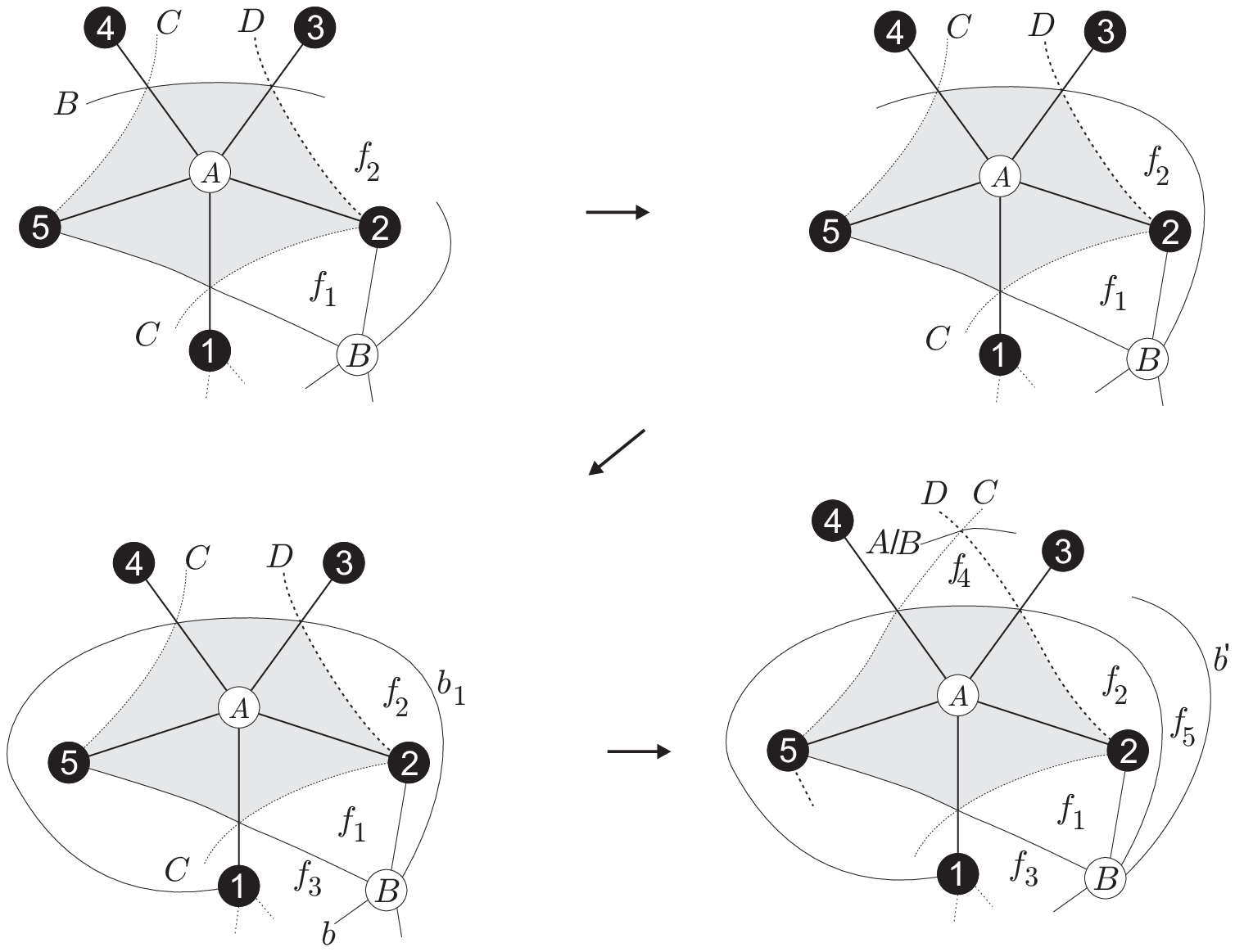}
\caption{$f_1$ is $3$-sided}
\label{fig:case1-1}
\end{figure}
Also, then $f_2$ cannot be $4$-sided, and thus $3$-sided. 
Hence the $B$-line, intersecting two $A$-lines $a_3$ and $a_4$, 
turns out to be $b_1$.
That is, it goes to black vertex $x_1$
as shown in the third of Figure \ref{fig:case1-1}.
(At this point, $b_1$ may contain triple crossing points on it after
crossing $a_4$.)

Consider the face $f_3$, which is adjacent to $f_1$ along the $B$-line $b_5$.
If $f_3$ is $3$-sided, then the $B$-line $b$, which must be $b_3$ or $b_4$, intersects the $A$-line $a_1$.
Then $b$ can reach neither $x_3$ nor $x_4$, because $b$ cannot
cross $b_1$ or meet $a_1$ twice.
Therefore, we found that $f_3$ is the only $4$-sided face.
Thus $f_4$ is $3$-sided.
We see that the line going through the upper triple crossing point of $f_4$
is either an $A$- or $B$-line by property ($\ast$).
Since $f_5$ is also $3$-sided, the $B$-line $b'$ goes to $x_3$, or crosses the $A$-line $a_3$.
In any case, $f_6$ cannot be $3$-sided
as described in Figure \ref{fig:case1-1-3}.
This is a contradiction.
\end{proof}

\begin{figure}[tb]
\includegraphics*[scale=0.7]{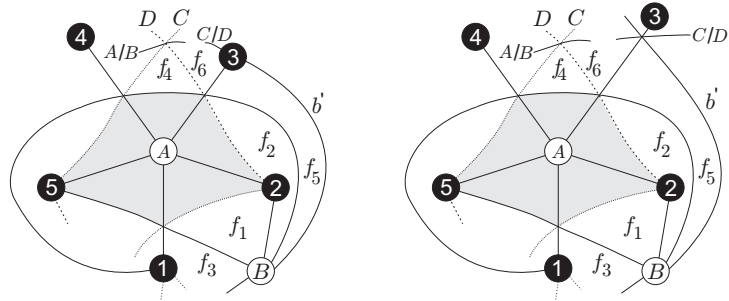}
\caption{$b'$ goes to $x_3$ or crosses $a_3$}
\label{fig:case1-1-3}
\end{figure}

\begin{claim}
Figure \ref{fig:case1-sub}(2) is impossible.
\end{claim}

\begin{proof}
First, since the horizontal $D$-line intersects $a_3$, $a_4$, $c_5$,
it is either $d_1$ or $d_2$ from our requirements for drawings.
By symmetry, we can assume that $f_1$ is $3$-sided as in Figure \ref{fig:case1-2}.
Then vertex $D$ is located around $f_1$.
Then the horizontal $D$-line turns out to be $d_1$.
However, neither $f_2$ nor $f_3$ is $3$-sided, which contradicts that
there is only one $4$-sided face.
\end{proof}

\begin{figure}[tb]
\includegraphics*[scale=0.7]{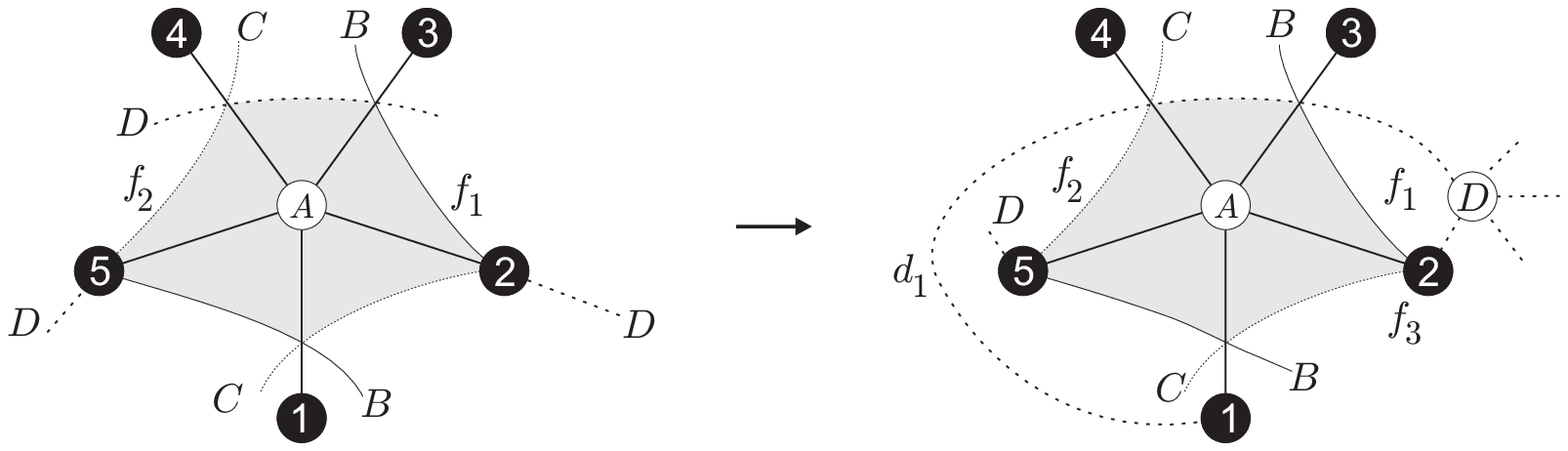}
\caption{$f_1$ is $3$-sided}
\label{fig:case1-2}
\end{figure}

\begin{claim}
Figure \ref{fig:case1-sub}(3) is impossible.
\end{claim}

\begin{proof}
In Figure \ref{fig:case1-3}, either $f_3$ or $f_4$ is  $4$-sided.
(For, if $f_3$ is not $4$-sided, then it is $3$-sided.
Then vertex $D$ is located there, which implies that $f_4$ is not $3$-sided.)
Thus both $f_1$ and $f_2$ are $3$-sided. 
Then the $B$-line $b_1$ is determined.
Since $f_5$ is $3$-sided, the $B$-line $b$, which is $b_3$ or $b_4$,
crosses the $A$-line $a_1$.
Then $b$ can reach neither $x_3$ nor $x_4$, a contradiction.
\end{proof}

\begin{figure}[tb]
\includegraphics*[scale=0.7]{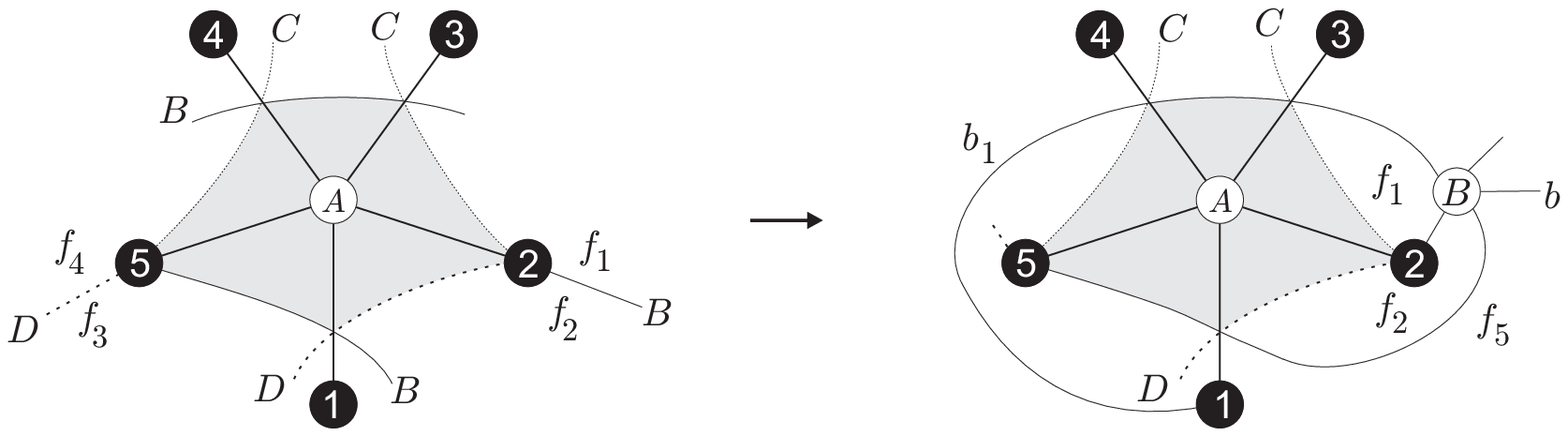}
\caption{Both $f_1$ and $f_2$ are $3$-sided}
\label{fig:case1-3}
\end{figure}

\begin{claim}
Figure \ref{fig:case1-sub}(4) is impossible.
\end{claim}

\begin{proof}
In Figure \ref{fig:case1-4},
if $f_1$ is not $3$-sided, then $f_2$ is $3$-sided, and
then vertex $D$ is located there.
But then, $f_1$ cannot be $4$-sided.
Hence $f_1$ is $3$-sided.
Similarly, so is $f_2$.
Then the $D$-line $d_1$ is determined.

\begin{figure}[tb]
\includegraphics*[scale=0.7]{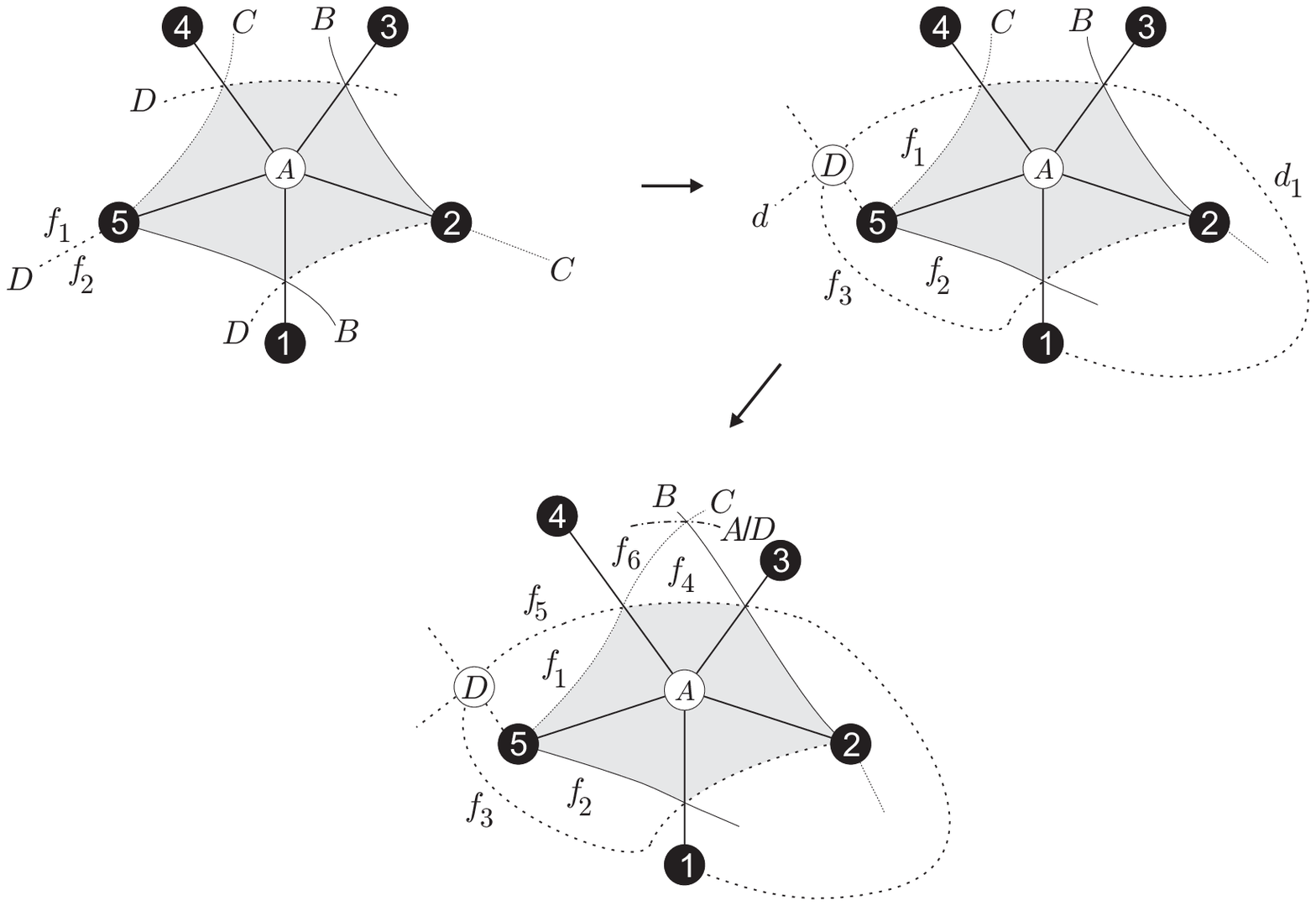}
\caption{Both $f_1$ and $f_2$ are $3$-sided}
\label{fig:case1-4}
\end{figure}

If $f_3$ is $3$-sided, then the $D$-line $d$, which is $d_3$ or $d_4$, crosses the $A$-line $a_1$.
Then $d$ cannot reach any black vertex as before.
Hence $f_3$ is $4$-sided.
As in the proof of Claim \ref{cl:case1}, examining $f_4,f_5,f_6$ leads to a contradiction.
\end{proof}

This completes the proof of Lemma \ref{lem:typeII-4}.
\end{proof}

\begin{lemma}\label{lem:typeII-2}
The number of type II triangles at vertex $A$ is not two.
\end{lemma}

\begin{proof}
Suppose that there are two type II triangles at $A$.
Then
we can assume that the local configuration at $A$ is as shown in Figure \ref{fig:case2start}(1),
up to renaming.
By property ($\ast$), the left upper line is a $B$- or $D$-line.
Similarly, the right upper line is a $C$- or $D$-line. 
See Figure \ref{fig:case2start}(2).
\begin{figure}[tb]
\includegraphics*[scale=0.7]{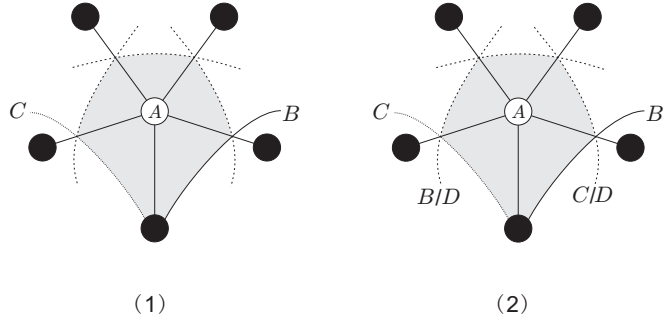}
\caption{Two type II triangles at $A$}
\label{fig:case2start}
\end{figure}
Then 
there are three cases, up to symmetry and relabeling of vertices, as shown in Figure \ref{fig:case2},
where the class of the horizontal line is determined by property ($\ast$).

\begin{figure}[tb]
\includegraphics*[scale=0.7]{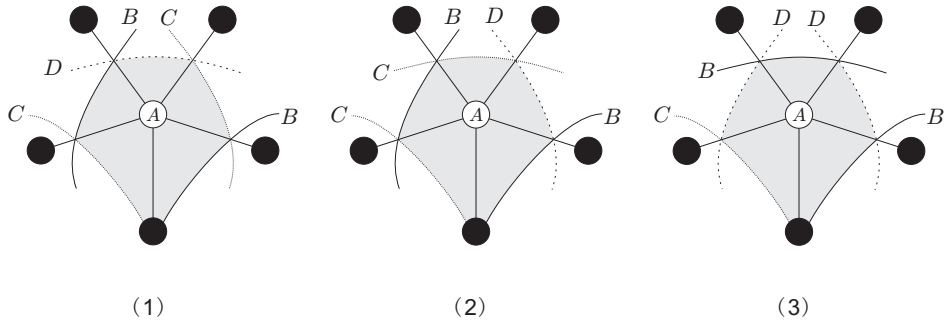}
\caption{Three cases where there are two type II triangles}
\label{fig:case2}
\end{figure}

\begin{claim}
Figure \ref{fig:case2}(1) is impossible.
\end{claim}

\begin{proof}
First, assume that $f_1$ is $4$-sided in Figure \ref{fig:case2-1-1}.
Then the others are all $3$-sided.
\begin{figure}[tb]
\includegraphics*[scale=0.7]{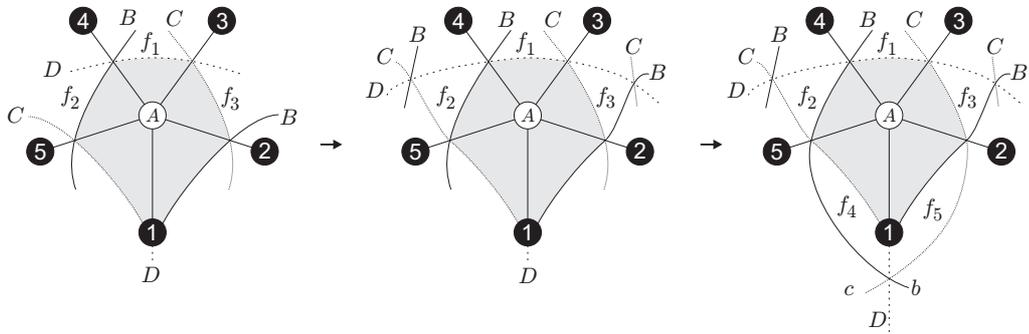}
\caption{Assume that $f_1$ is $4$-sided}
\label{fig:case2-1-1}
\end{figure}
Thus $f_2$ and $f_3$, and then $f_4$, $f_5$ are determined as in Figure \ref{fig:case2-1-1}.
(If an $A$-line goes through the left triple crossing point of $f_2$, then
the face sharing a $D$-line with $f_2$ cannot be $3$-sided. Similarly for $f_3$.)

Consider the $B$-line $b$.
It goes to $x_2$ or crosses the $A$-line $a_2$.
Suppose that the former happens.
Then $f_6,\dots,f_9$ are determined as in Figure \ref{fig:case2-1-1a}.
Moreover, the $D$-line $d_5$ is also determined.

\begin{figure}[tb]
\includegraphics*[scale=0.7]{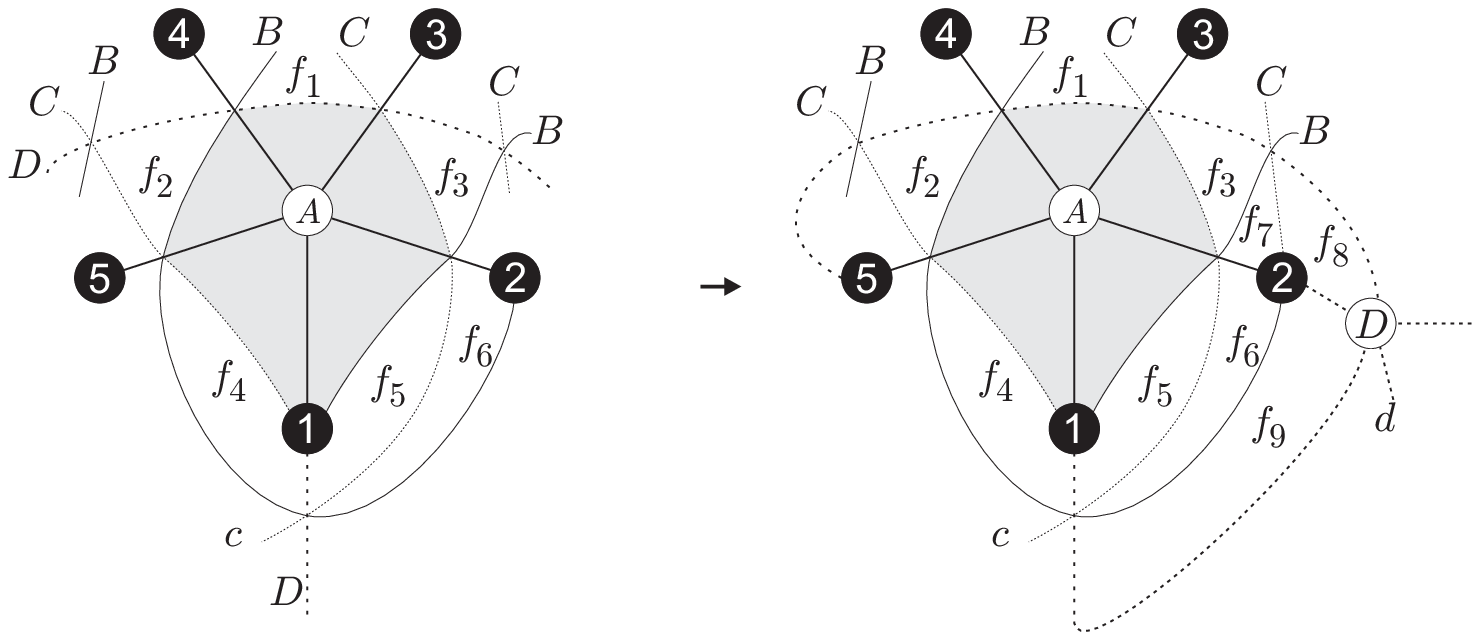}
\caption{The case where $b$ goes to $x_2$}
\label{fig:case2-1-1a}
\end{figure}

Then the $C$-line $c$ cannot go to $x_5$, since it crosses $d_5$.
Hence it crosses the $A$-line $a_5$.
This forces the $D$-line $d$ to cross the same $a_5$.
Then it cannot reach any black vertex, a contradiction.
Therefore, $b$ crosses the $A$-line $a_2$.
By the same reason, $c$ crosses $a_5$.

Repeating the same argument, we obtain the configuration as shown in
Figure \ref{fig:case2-1-1b}. 
\begin{figure}[tb]
\includegraphics*[scale=0.7]{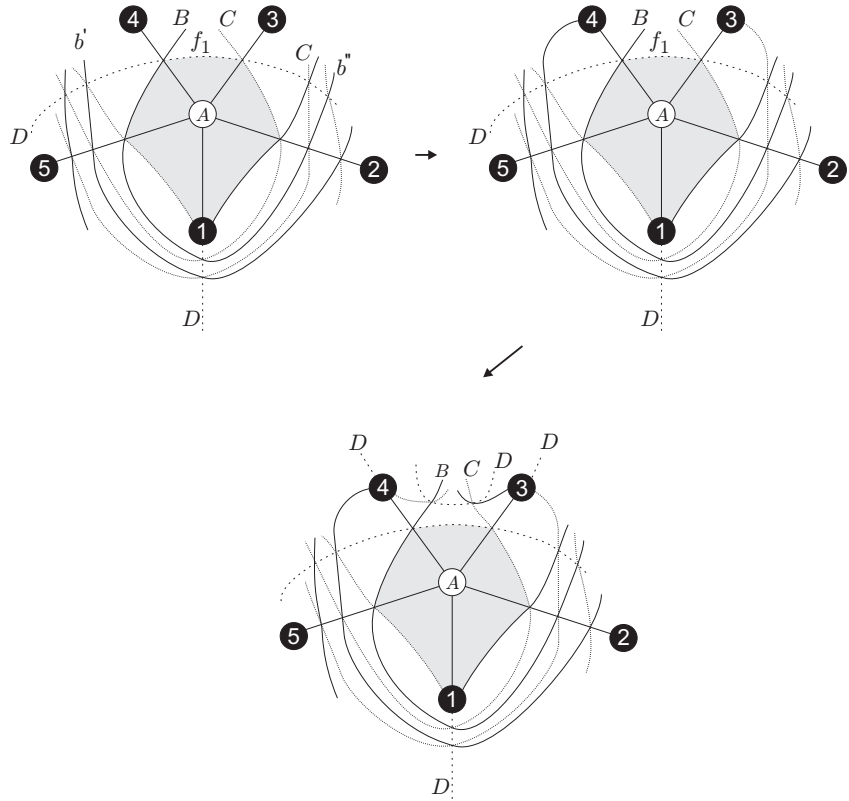}
\caption{A final contradiction when $f_1$ is $4$-sided}
\label{fig:case2-1-1b}
\end{figure}
If the $B$-line $b'$ crosses the $A$-line $a_4$,
then both $b'$ and $b''$ go to $x_3$, a contradiction.
Thus $b'$ goes to $x_4$.
Similarly, the $C$-line $c_3$ is determined.
See the second of Figure \ref{fig:case2-1-1b}.

Then $f_1$ can be incident with neither vertex $B$ nor $C$.
For example, if $f_1$ is incident with $B$, then
the left face of $f_1$ cannot be $3$-sided.
Hence $f_1$ is incident with two more triple crossing points.
Then the upper horizontal line of $f_1$ is an $A$- or $D$-line by property ($\ast$).
From our requirements, it cannot be an $A$-line.
Thus we have the third of Figure \ref{fig:case2-1-1b},
but then vertex $D$ cannot be located.

Next, assume that $f_1$ is $3$-sided.
We see that a $D$-line goes through
the upper triple crossing point of $f_1$ from our requirements and property ($\ast$).

By symmetry, we can assume that $f_2$ is $3$-sided.
See Figure \ref{fig:case2-1-2}.
\begin{figure}[tb]
\includegraphics*[scale=0.7]{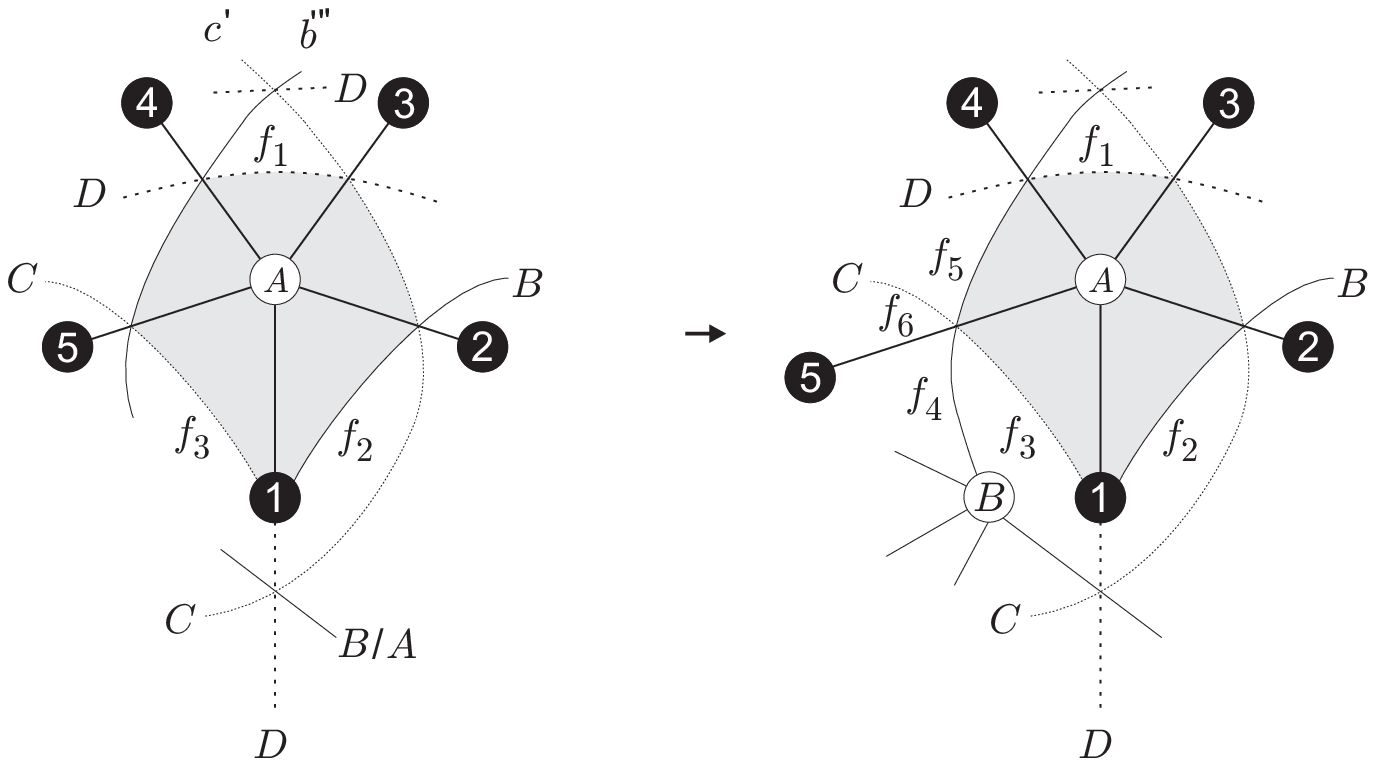}
\caption{$f_1$ is $3$-sided}
\label{fig:case2-1-2}
\end{figure}
If $f_3$ is also $3$-sided, then the $B$-line $b'''$ and the $C$-line $c'$ meet
twice, a contradiction.
Hence $f_3$ turns out to be $4$-sided.
Also, no $A$-line is adjacent to $f_3$, because
each of the four $A$-lines $a_2,\dots,a_5$ meets $b'''$ or $c'$. 
Thus vertex $B$ is located as in Figure \ref{fig:case2-1-2}.
Again, examining $f_4,f_5,f_6$ leads to a contradiction as in the proof of Claim \ref{cl:case1}.
\end{proof}

\begin{claim}
Figure \ref{fig:case2}(2) is impossible.
\end{claim}

\begin{proof}
In Figure \ref{fig:case2-2},
suppose that $f_1$ is $4$-sided.
Then $f_2$ is $3$-sided, and so vertex $D$ appears there.
Then $f_3$ is not $3$-sided, a contradiction.
Hence $f_1$ is $3$-sided.

\begin{figure}[tb]
\includegraphics*[scale=0.7]{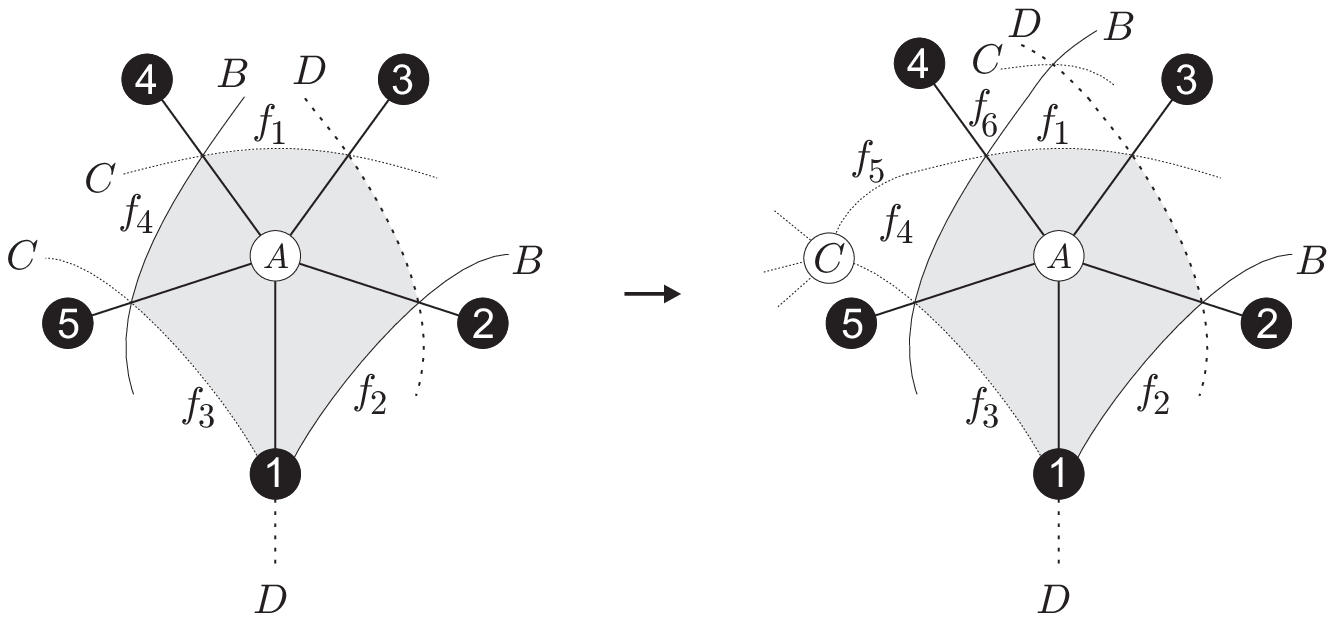}
\caption{$f_1$ is $3$-sided}
\label{fig:case2-2}
\end{figure}

If $f_2$ is $3$-sided, then $f_3$ is $4$-sided as above.
Otherwise, $f_2$ is $4$-sided.
In any case, $f_4$ is $3$-sided, and vertex $C$ appears.
Also, $f_5$ and $f_6$ are $3$-sided.
But this is impossible as in the proof of Claim \ref{cl:case1} again.
\end{proof}

\begin{claim}
Figure \ref{fig:case2}(3) is impossible.
\end{claim}

\begin{proof}
In Figure \ref{fig:case2-3}, at least two of $f_1, f_2, f_3$ are $3$-sided.
If $f_1$ and $f_2$ are $3$-sided, then vertex $D$ cannot be located correctly.
Similarly for the case where $f_1$ and $f_3$ are $3$-sided.
Hence $f_2$ and $f_3$ are $3$-sided.
Then the $D$-line $d$ goes to $x_2$ or $x_3$, and another $D$-line $d'$ goes to
$x_4$ or $x_5$.
But this is impossible.
\end{proof}

\begin{figure}[tb]
\includegraphics*[scale=0.7]{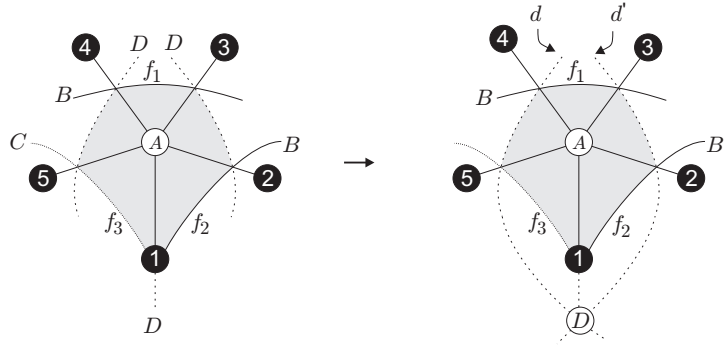}
\caption{$f_2$ and $f_3$ are $3$-sided}
\label{fig:case2-3}
\end{figure}

This completes the proof of Lemma \ref{lem:typeII-2}.
\end{proof}

\begin{lemma}\label{lem:typeII-0}
The number of type II triangles at vertex $A$ is not zero.
\end{lemma}

\begin{proof}
Assume that there are no type II triangles at $A$.
Up to symmetry and relabeling of vertices,
the local configuration at $A$ can be assumed as in the first of Figure \ref{fig:case3}.

\begin{figure}[tb]
\includegraphics*[scale=0.7]{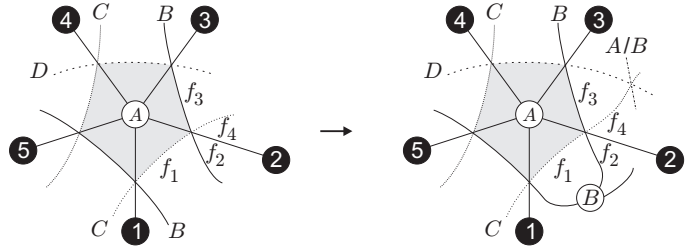}
\caption{Five type I triangles at $A$}
\label{fig:case3}
\end{figure}

By symmetry, we can assume that the right hand side does not contain a $4$-sided face.
More precisely, $f_1,\dots, f_4$ are all $3$-sided.
Thus vertex $B$ is located.
Then examining $f_2,f_3,f_4$, as in the proof of Claim \ref{cl:case1}, leads to a contradiction.
(In this case, $f_3$ can be incident with $x_5$.
Then $f_4$ cannot be $3$-sided likewise.)
\end{proof}

\begin{theorem}\label{thm:k45}
$K_{5,4}$ does not admit a semi-regular drawing.
\end{theorem}

\begin{proof}
This follows from Lemmas 
\ref{lem:typeII-4}, \ref{lem:typeII-2} and \ref{lem:typeII-0}.
\end{proof}

\section{$K_{4,4}$}\label{sec:k44}

Throughout this section, we will assume that $G=K_{4,4}$ admits a semi-regular drawing.
We will show that this is impossible.

Let 
$V_1=\{x_1,x_2,x_3,x_4\}$ and $V_2=\{A,B,C,D\}$ be the partite sets of $G$.
As in Section \ref{sec:k54}, we refer to vertices of $V_1$ (resp. $V_2$) as black (resp. white) vertices, and use the same notions as
$A$-lines, type I or II triangles, so on.
Also, property ($\ast$) holds from our requirements.

We fix a semi-regular drawing of $G$, which is denoted by $G$ again.
Let $k$ be the number of triple crossing points.
Add a new vertex to each triple crossing point.
Then we have a plane graph $G'$ with $8+k$ vertices and $16+3k$ edges.
Since $3(8+k)-(16+3k)-6=2$,
either
\begin{enumerate}
\item one face of $G'$ is $5$-sided, and the others are $3$-sided; or
\item two faces of $G'$ are $4$-sided, and the others are $3$-sided
\end{enumerate}
by Lemma \ref{lem:faces}.
As in Section \ref{sec:k54},
a face of $G$ means that of $G'$.

\subsection{Case (1)}
We treat the case where one face of $G$ is $5$-sided, and the others are $3$-sided.
At most two white vertices appear in the $5$-sided face.
Hence we can assume that four faces at vertex $A$ are all $3$-sided.
Thus the number of type II triangles at vertex $A$ is either $0$, $2$ or $4$.
We will eliminate these three possibilities.

\begin{lemma}\label{lem:i4}
The number of type II triangles at vertex $A$ is not four.
\end{lemma}

\begin{proof}
Assume that there are four type II triangles at $A$.
We may assume that the local configuration at $A$ is as shown in Figure \ref{fig:casei-0}(1),
up to renaming.
By property ($\ast$), the right upper line is a $C$- or $D$-line,
and the left upper line is a $B$- or $D$-line (see Figure \ref{fig:casei-0}(2)).
\begin{figure}[tb]
\includegraphics*[scale=0.7]{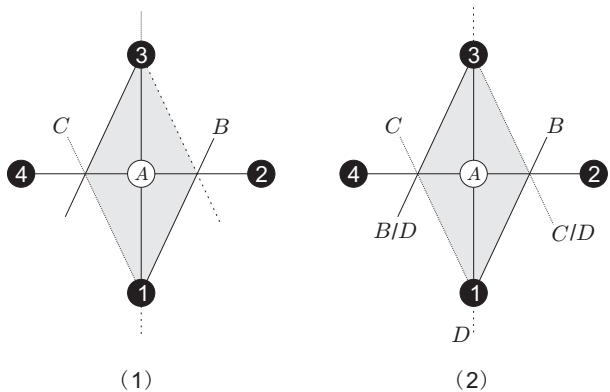}
\caption{Four type II triangles at $A$}
\label{fig:casei-0}
\end{figure}
Up to symmetry and relabeling, there are two possibilities as shown in Figure \ref{fig:casei-1}.

\begin{figure}[tb]
\includegraphics*[scale=0.7]{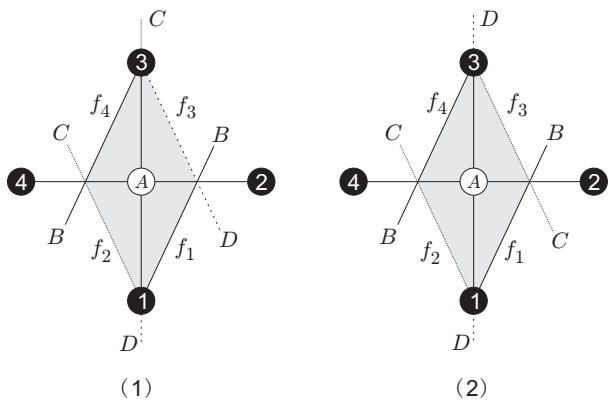}
\caption{Two cases where there are two type II triangles}
\label{fig:casei-1}
\end{figure}

\begin{claim}
Figure \ref{fig:casei-1}(1) is impossible.
\end{claim}

\begin{proof}
By symmetry, we can assume that both $f_1$ and $f_2$ are $3$-sided.
Since $f_1$ is $3$-sided, vertex $D$ is located.
Then $f_2$ cannot be $3$-sided, a contradiction.
\end{proof}

\begin{claim}\label{cl:cc}
Figure \ref{fig:casei-1}(2) is impossible.
\end{claim}

\begin{proof}
By symmetry, we can assume that $f_1$, $f_2$ are $3$-sided again.
Then the $B$-line $b_3$ meets the $C$-line $c_3$, a contradiction.
\end{proof}

This complete the proof of Lemma \ref{lem:i4}.
\end{proof}

\begin{lemma}\label{lem:i2}
The number of type II triangles at vertex $A$ is not two.
\end{lemma}

\begin{proof}
Assume that there are two type II triangles at $A$.
As before, we may assume that the local configuration at $A$ is as shown in
Figure \ref{fig:casei-2start}(1).
By property ($\ast$), the right upper line is a $C$- or $D$-line,
and the left upper line is a $B$- or $D$-line (see Figure \ref{fig:casei-2start}).
Up to symmetry and relabeling, there are two possibilities as shown in Figure \ref{fig:casei-2}.

\begin{figure}[tb]
\includegraphics*[scale=0.7]{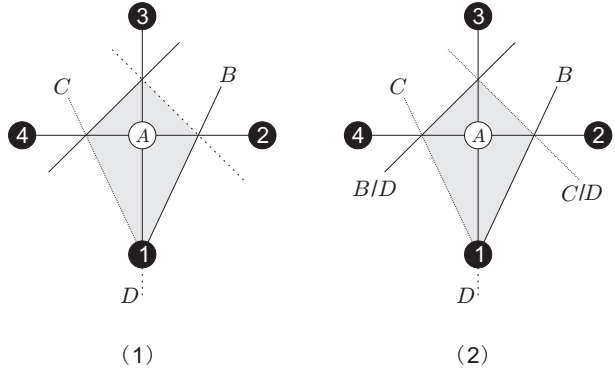}
\caption{Two type II triangles at $A$}
\label{fig:casei-2start}
\end{figure}

\begin{figure}[tb]
\includegraphics*[scale=0.7]{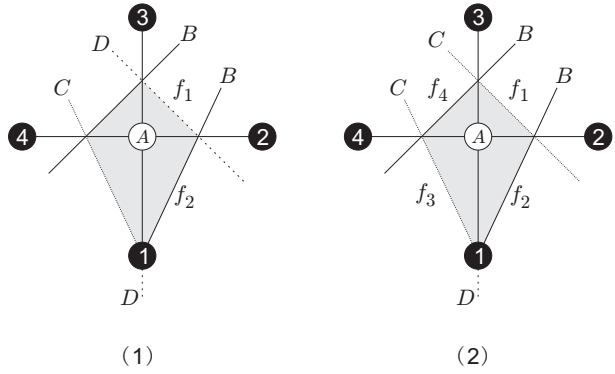}
\caption{Two cases where there are two type II triangles}
\label{fig:casei-2}
\end{figure}

\begin{claim}
Figure \ref{fig:casei-2}(1) is impossible.
\end{claim}

\begin{proof}
Assume that $f_2$ is not $3$-sided.
Then $f_1$ is $3$-sided, and so vertex $B$ is located there.
Thus the $B$-line $b_2$ is determined.
See the first of Figure \ref{fig:casei-2-1}.
Another $B$-line $b$ crosses the $A$-line $a_2$, but
then it cannot reach any black vertex.

\begin{figure}[tb]
\includegraphics*[scale=0.7]{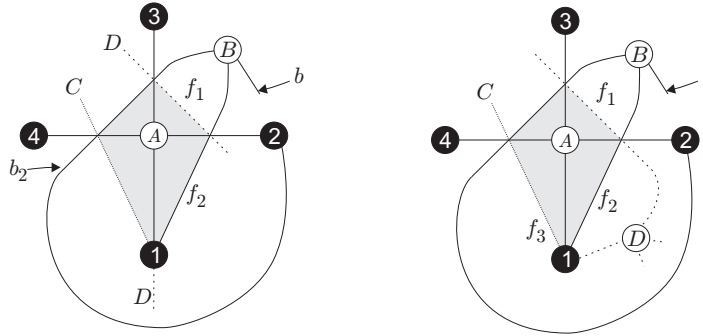}
\caption{$b$ crosses $a_2$}
\label{fig:casei-2-1}
\end{figure}

Hence $f_2$ is $3$-sided, so vertex $D$ is located.
Then $f_3$ cannot be $3$-sided. 
Thus $f_1$ is $3$-sided, so vertex $B$ is located, and the $B$-line $b_2$ is determined again
as in the second of Figure \ref{fig:casei-2-1}.
Examining $b$ leads to a contradiction as above.
\end{proof}

\begin{claim}\label{cl:i-2}
Figure \ref{fig:casei-2}(2) is impossible.
\end{claim}

\begin{proof}
By symmetry, we may assume that both $f_1$ and $f_2$ are $3$-sided.
Then vertex $B$ is located.  See the first of Figure \ref{fig:casei-2-2}.

\begin{figure}[tb]
\includegraphics*[scale=0.7]{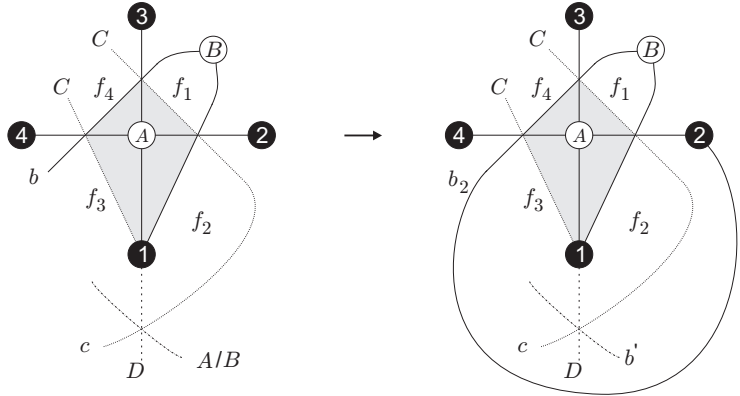}
\caption{Both $f_1$ and $f_2$ are $3$-sided}
\label{fig:casei-2-2}
\end{figure}

If $f_3$ is $3$-sided, then 
the $B$-line $b$ and the $C$-line $c$ meet twice, a contradiction.
Hence $f_3$ is not $3$-sided, and $b$ goes to $x_2$ as in Figure \ref{fig:casei-2-2}.
Then a line $b'$ turns out to be a $B$-line.
But this $B$-line cannot reach any black vertex, otherwise it crosses $a_2$ twice.
\end{proof}

This completes the proof of Lemma \ref{lem:i2}.
\end{proof}

\begin{lemma}\label{lem:i0}
The number of type II triangles at vertex $A$ is not zero.
\end{lemma}

\begin{proof}
Assume that there is no type II triangle at $A$.
Up to symmetry and relabeling, there are two possibilities as shown in Figure \ref{fig:casei-3}.

\begin{figure}[tb]
\includegraphics*[scale=0.7]{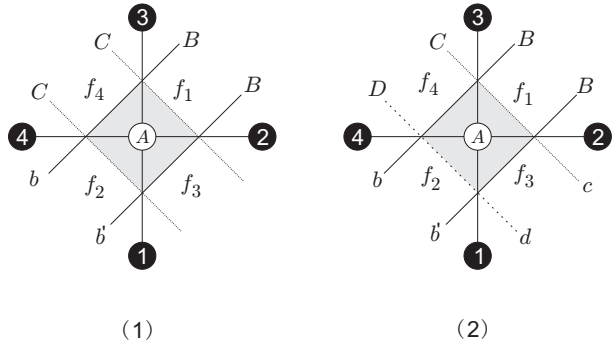}
\caption{Two cases where there are four type I triangles}
\label{fig:casei-3}
\end{figure}

In any case, we can assume that $f_1$ is $3$-sided by symmetry.
Then vertex $B$ is located there.
The $B$-line $b$ goes to $x_1$ or $x_2$, and
the $B$-line $b'$ goes to $x_3$ or $x_4$.
This is impossible.
\end{proof}

\begin{proposition}\label{prop:casei} 
Case (1) is impossible.
\end{proposition}

\begin{proof}
This immediately follows from Lemmas \ref{lem:i4}, \ref{lem:i2} and \ref{lem:i0}.
\end{proof}

\subsection{Case (2)}

We treat the case where
two faces of $G$ are $4$-sided, and the others are $3$-sided.
There are two subcases.

\begin{itemize}
\item[(2-1)] All white vertices are incident with a $4$-sided face.
\item[(2-2)] There is a white vertex which is not incident with a $4$-sided face.
\end{itemize}

\subsubsection{Subcase (2-1).}

Around each $4$-sided face, just two white vertices appear.
Hence there is just one $4$-sided face around each white vertex.
Also, it implies that
if a face is incident with two adjacent triple crossing points or
an adjacent pair of a triple crossing point and a black vertex,
then it must be $3$-sided.

In this subcase, then the number of type II triangle at vertex $A$ is $0$, $1$, $2$ or $3$.

\begin{lemma}\label{lem:ii3}
The number of type II triangles at $A$ is not three.
\end{lemma}

\begin{proof}
Assume that there are thee type II triangle at $A$.
Up to symmetry and renaming, 
the situation is as shown in Figure \ref{fig:caseii3}(1),
where $f$ is $4$-sided.

\begin{figure}[tb]
\includegraphics*[scale=0.7]{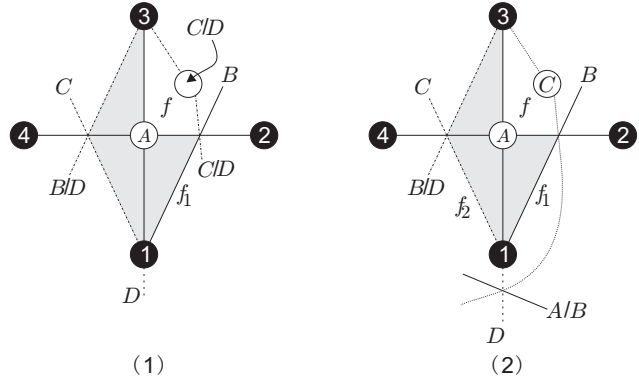}
\caption{Three type II triangles at $A$}
\label{fig:caseii3}
\end{figure}

As remarked above, $f_1$ is $3$-sided.
If $f$ is incident with vertex $D$,
then $f_1$ is incident with vertex $D$.
This is impossible.
Thus $f$ is incident with vertex $C$.
See Figure \ref{fig:caseii3}(2).
Since $f_2$ is also $3$-sided, 
the left upper line of $f_2$ is a $B$-line.
We have the configuration as in Figure \ref{fig:caseii3a}(1).
\begin{figure}[tb]
\includegraphics*[scale=0.7]{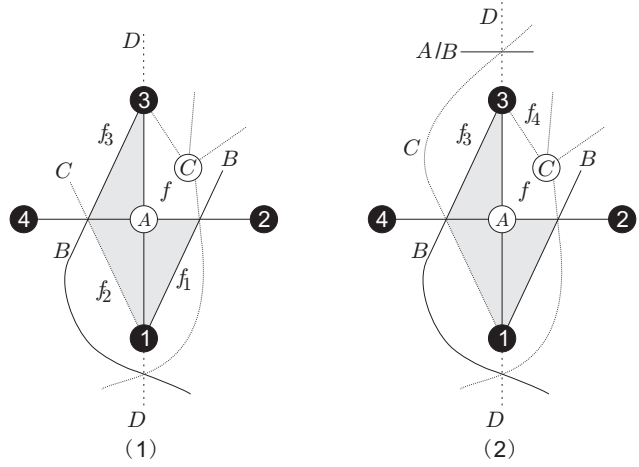}
\caption{Three type II triangles at $A$ (continued)}
\label{fig:caseii3a}
\end{figure}
Then $f_3$ is $3$-sided, but this forces 
$f_4$ to be neither $3$-sided nor $4$-sided, a contradiction.
See Figure \ref{fig:caseii3a}(2). 
\end{proof}

\begin{lemma}\label{lem:ii2}
The number of type II triangles at $A$ is not two.
\end{lemma}

\begin{proof}
Assume that there are two type II triangles at $A$.
Up to symmetry and renaming, there are two possibilities as shown in Figure \ref{fig:caseii-1-1},
where $f$ is $4$-sided.
Notice that $f_1$ and $f_2$ are $3$-sided.

\begin{figure}[tb]
\includegraphics*[scale=0.7]{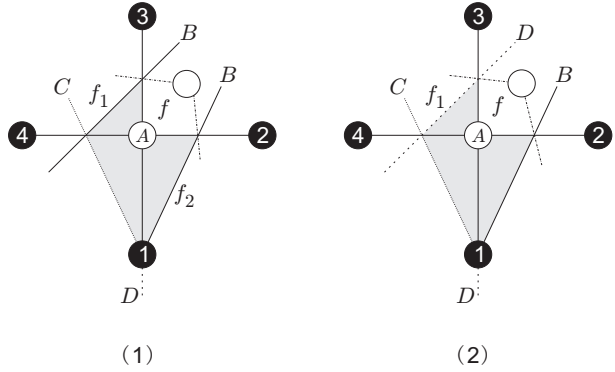}
\caption{Two cases where there are two type II triangles}
\label{fig:caseii-1-1}
\end{figure}

For Figure \ref{fig:caseii-1-1}(1),
either vertex $C$ or $D$ appears around $f$ by property ($\ast$).
Since $f_1$ is $3$-sided, the former is impossible.
The latter is also impossible, because $f_2$ is $3$-sided.
For Figure \ref{fig:caseii-1-1}(2),
vertex $C$ appears around $f$ by property ($\ast$).
Then $f_1$ gives a contradiction, again.
\end{proof}

\begin{lemma}\label{lem:ii1}
The number of type II triangles at $A$ is not one.
\end{lemma}

\begin{proof}
Suppose that there is one type II triangle at $A$.
Up to symmetry and renaming, the configuration is as shown in Figure \ref{fig:caseii1},
where $f$ is $4$-sided.

\begin{figure}[tb]
\includegraphics*[scale=0.7]{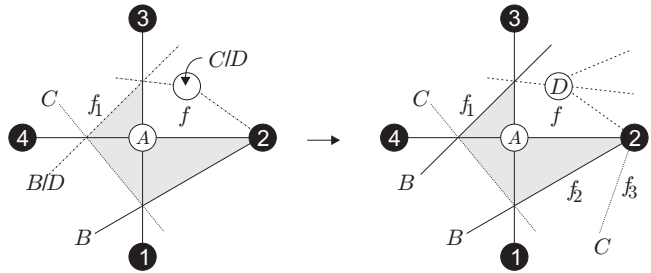}
\caption{One type II triangle at $A$}
\label{fig:caseii1}
\end{figure}

Notice that $f_1$ is $3$-sided.
If $f$  is incident with vertex $C$, then
$f_1$ is incident with $C$, an impossible.
Hence $f$ is incident with vertex $D$.
Thus we have the configuration as in Figure \ref{fig:caseii1}.
The fact that $f_2$ is $3$-sided forces $f_3$ to be $4$-sided.
Then $f_3$ is incident with 
vertices $C$ and $D$, which
contradicts the fact that $D$ is incident with only one $4$-sided face.
\end{proof}

\begin{lemma}\label{lem:ii0}
The number of type II triangles at $A$ is not zero.
\end{lemma}

\begin{proof}
Assume that there is no type II triangle at $A$.
Up to symmetry and renaming, there are two possibilities as shown in Figure \ref{fig:caseii-1-2},
where $f$ is $4$-sided.

\begin{figure}[tb]
\includegraphics*[scale=0.7]{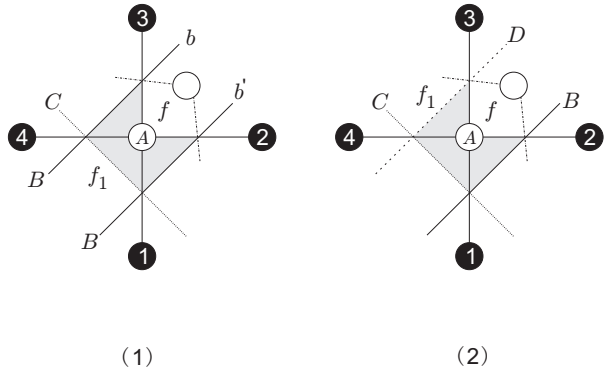}
\caption{Two cases where there is no type I triangle}
\label{fig:caseii-1-2}
\end{figure}

For Figure \ref{fig:caseii-1-2}(1),
$f_1$ is $3$-sided, so vertex $B$ is located there.
Then $b$ goes to $x_1$ or $x_2$, and $b'$ goes to $x_3$ or $x_4$.
This is impossible.

For Figure \ref{fig:caseii-1-2}(2),
vertex $C$ appears around $f$ by property ($\ast$).
On the other hand, $f_1$ is $3$-sided.
This is impossible.
\end{proof}

\subsubsection{Subcase (2-2).}

In this subcase, we may assume that vertex $A$ is not incident with a $4$-sided face
without loss of generality.
Thus the number of type II triangles at $A$ is $0$, $2$ or $4$.

\begin{lemma}\label{lem:ii-2-4}
The number of type II triangles at $A$ is not four.
\end{lemma}

\begin{proof}
Suppose that there are four type II triangles at $A$.
Then
there are two possibilities as in Figure \ref{fig:casei-1}.

\begin{claim}
Figure \ref{fig:casei-1}(1) is impossible.
\end{claim}

\begin{proof}
Among the four faces $f_1,\dots, f_4$,
at least two are $3$-sided.
Furthermore, if $f_1$ (resp. $f_3$) is $3$-sided, then $f_2$ (resp. $f_4$) is $4$-sided, and vice versa.
Up to symmetry,
there are three possibilities:
\begin{itemize}
\item[(a)] $f_1$ and $f_3$ are $3$-sided.
\item[(b)] $f_1$ and $f_4$ are $3$-sided.
\item[(c)] $f_2$ and $f_3$ are $3$-sided.
\end{itemize}

(a) In this case, $f_2$ and $f_4$ are $4$-sided.
Thus the others are all $3$-sided.
See Figure \ref{fig:caseii-2-1-a}(1).
By examining $f_5$, the left upper line of $f_4$ is not an $A$-line.
Since $f_5$ and $f_6$ are $3$-sided, two $D$-lines $d$ and $d'$ go to $x_2$, or
cross, a contradiction.

\begin{figure}[tb]
\includegraphics*[scale=0.7]{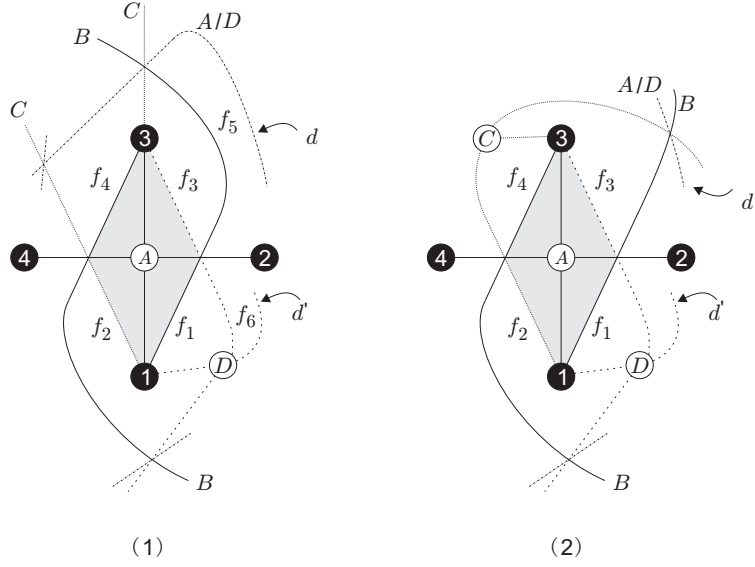}
\caption{For (a) and (b)}
\label{fig:caseii-2-1-a}
\end{figure}

(b) In this case, $f_2$ and $f_3$ are $4$-sided.
By examining $d$ and $d'$ shown in Figure \ref{fig:caseii-2-1-a}(2),
the same argument as (a) leads to a contradiction.

(c) In this case, $f_1$ and $f_4$ are $4$-sided.
As above, we see that neither $f_1$ nor $f_4$ is incident with an $A$-line.
See the first of Figure \ref{fig:caseii-2-1-c}.
\begin{figure}[tb]
\includegraphics*[scale=0.7]{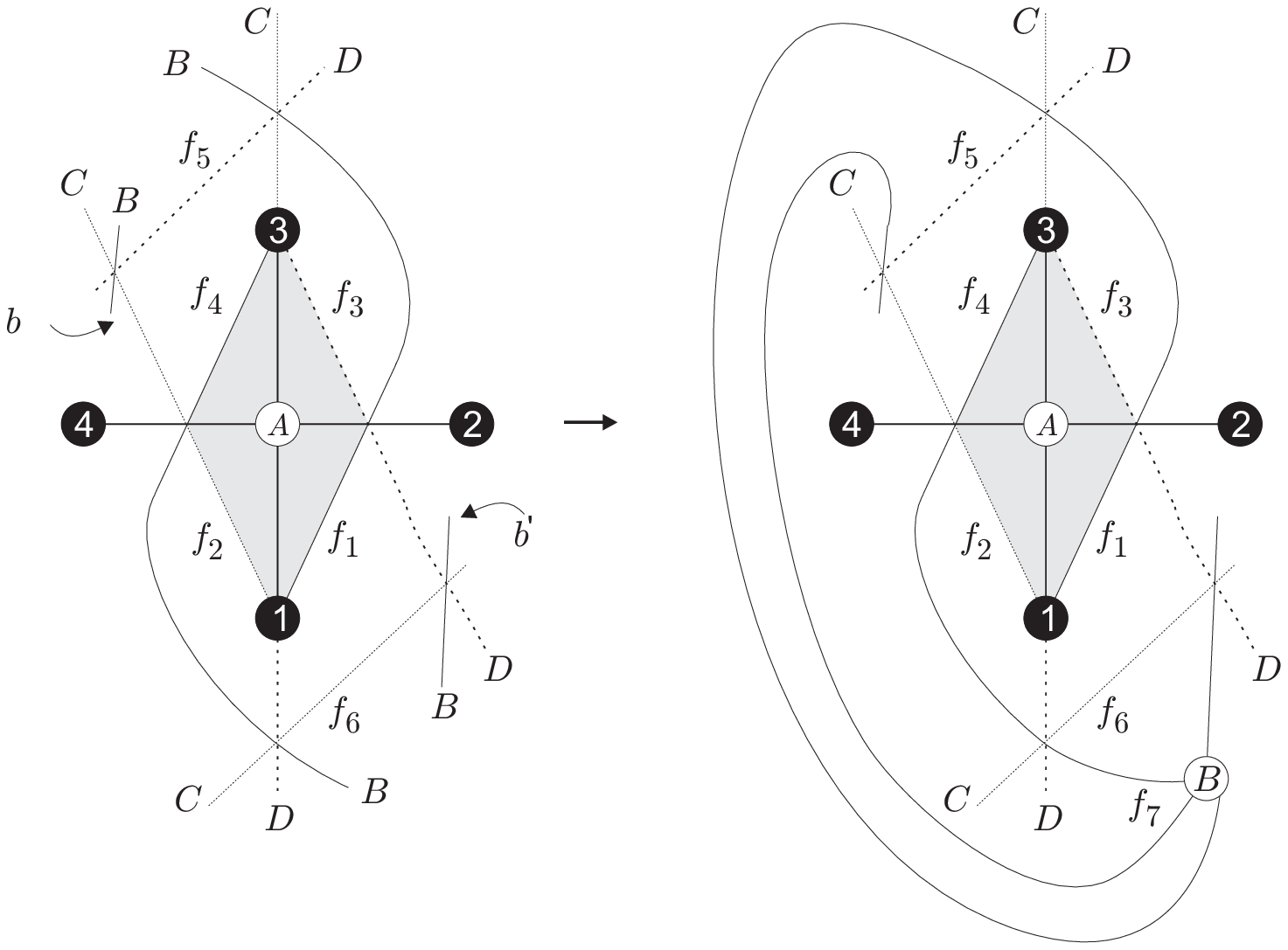}
\caption{For (c)}
\label{fig:caseii-2-1-c}
\end{figure}
Thus we have two $B$-lines $b$ and $b'$ as shown there.
Since $f_5$ and $f_6$ are $3$-sided, vertex $B$ is located as in the second of Figure \ref{fig:caseii-2-1-c}.
But then $f_7$ cannot be $3$-sided, a contradiction.
\end{proof}

\begin{claim}
Figure \ref{fig:casei-1}(2) is impossible.
\end{claim}

\begin{proof}
If both $f_1$ and $f_2$ are $3$-sided, then
we have a contradiction as in the proof of Claim \ref{cl:cc}.
Hence either of $f_1$ or $f_2$ is $4$-sided.
Similarly, either $f_3$ or $f_4$ is $4$-sided.
Then there are two possibilities, up to symmetry.

\begin{itemize}
\item[(d)] $f_1$ and $f_3$ are $4$-sided.
\item[(e)] $f_1$ and $f_4$ are $4$-sided.
\end{itemize}

(d) Then $f_2$ and $f_4$ are $3$-sided as in the first of Figure \ref{fig:caseii-2-1-1}.
\begin{figure}[tb]
\includegraphics*[scale=0.65]{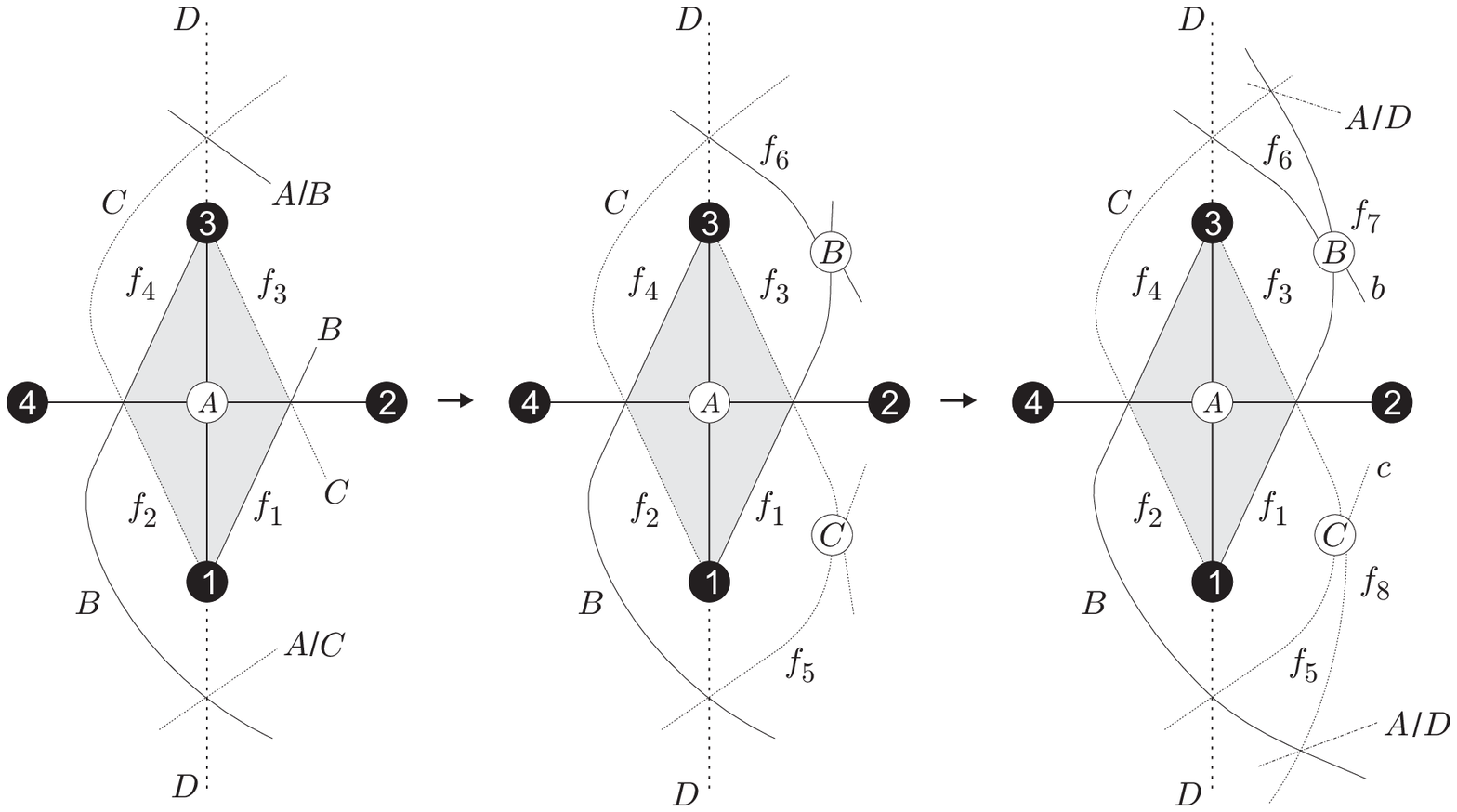}
\caption{For (d)}
\label{fig:caseii-2-1-1}
\end{figure}
By property ($\ast$), the right upper line $\ell$ of $f_3$
is an $A$- or $B$-line.
If $\ell$ is an $A$-line, then it is either $a_2$ or $a_4$.
But if $\ell$ is $a_4$, then $\ell$ meets $c_1$ twice, impossible.
If $\ell$ is $a_2$, then $\ell$ meets $b_1$ twice, impossible.
Thus $\ell$ is a $B$-line.
By the same reason, the right lower line of $f_1$ is not an $A$-line, and 
so a $C$-line.

Thus vertices $B$ and $C$ are located as in the second of Figure \ref{fig:caseii-2-1-1}.
After locating $f_5$ and $f_6$ as in the third of Figure \ref{fig:caseii-2-1-1},
consider the $B$-line $b$ and the $C$-line $c$. 
If $b$ crosses the $A$-line $a_2$, then so does $c_1$ through the same triple crossing point on $a_2$.
Then $f_7$ cannot be $3$-sided.
Hence $b$, and then $c$,  go to $x_2$.
Then $f_7$ and $f_8$ cannot be $3$-sided simultaneously.

(e) Then $f_2$ and $f_3$ are $3$-sided.
See Figure \ref{fig:caseii-2-1-2}.
\begin{figure}[tb]
\includegraphics*[scale=0.7]{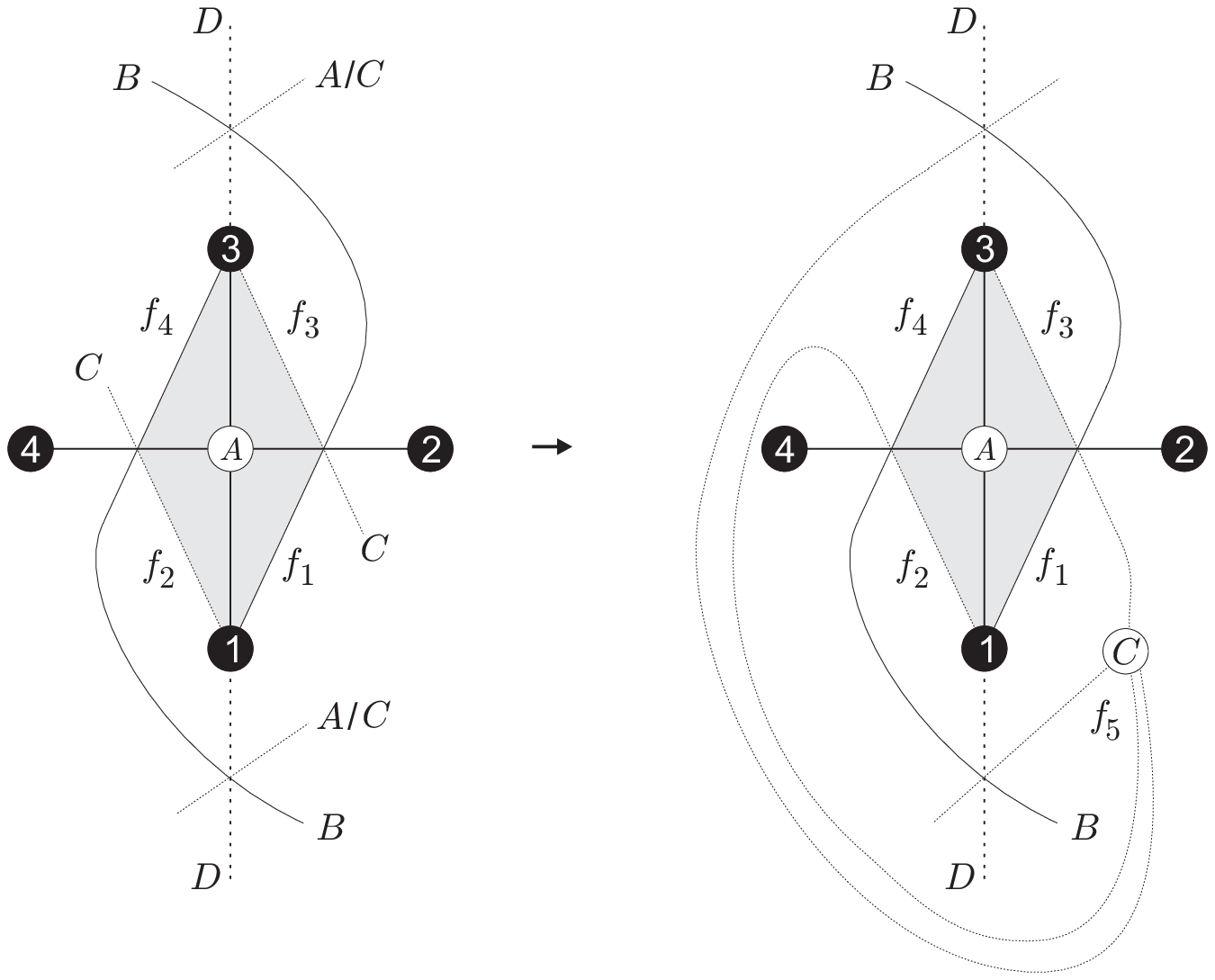}
\caption{For (e)}
\label{fig:caseii-2-1-2}
\end{figure}

By the same argument as (d), 
the right lower line of $f_1$ turns out to be a $C$-line.
Similarly, the upper line of $f_4$ is a $C$-line.
Since both $f_1$ and $f_4$ are $4$-sided, 
vertex $C$ is located as in the second of Figure \ref{fig:caseii-2-1-2}.
Then $f_5$ is not $3$-sided, a contradiction.
\end{proof}

This completes the proof of Lemma \ref{lem:ii-2-4}.
\end{proof}

\begin{lemma}\label{lem:ii-2-2}
The number of type II triangles at $A$ is not two.
\end{lemma}

\begin{proof}
Suppose that there are two type II triangles at $A$.
Then
the local configuration at $A$ is Figure \ref{fig:casei-2}(1) or (2).

\begin{claim}\label{cl:ii-2-2}
Figure \ref{fig:casei-2}(1) is impossible.
\end{claim}

\begin{proof}
We claim that $f_1$ is $4$-sided.
Assume that $f_1$ is $3$-sided.
Then vertex $B$ is located, and the $B$-line $b_2$ is
determined as in Figure \ref{fig:caseii-2-2-1}.

\begin{figure}[tb]
\includegraphics*[scale=0.7]{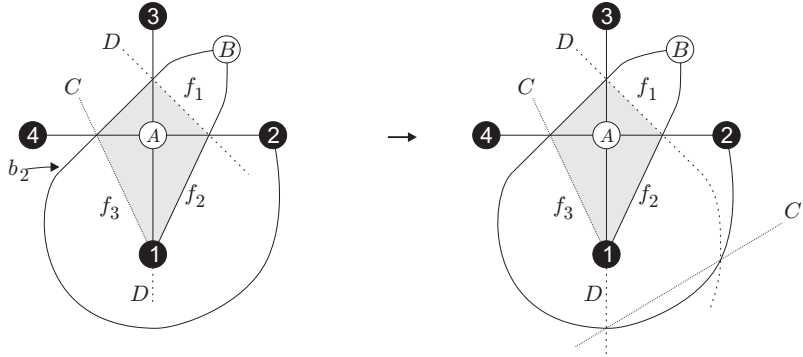}
\caption{The case where $f_1$ is $3$-sided}
\label{fig:caseii-2-2-1}
\end{figure}

Suppose further that $f_2$ is $4$-sided.
If vertex $D$ is incident with $f_2$,
then $f_2$ cannot be $4$-sided.
Hence $f_2$ is incident with two more triple crossing points.
Then the fourth line of $f_2$ is an $A$-, $B$- or $C$-line by property ($\ast$).
However, the existence of $b_2$ implies that it is neither an $A$- nor $B$-line.
Thus the right lower line of $f_2$ is a $C$-line.
Then the second of Figure \ref{fig:caseii-2-2-1} is the only possible configuration for $f_2$.
But this is impossible, because two lines meet at most once.

\begin{figure}[tb]
\includegraphics*[scale=0.7]{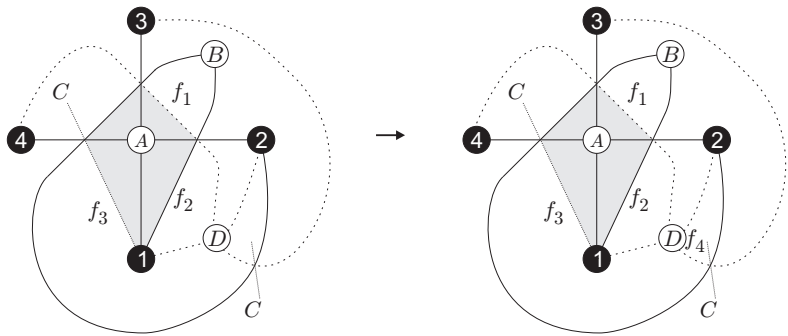}
\caption{The case where $f_1$ is $3$-sided (continued)}
\label{fig:caseii-2-2-1-1}
\end{figure}

Thus we see that $f_2$ is $3$-sided, so vertex $D$ is located as in Figure \ref{fig:caseii-2-2-1-1}.
Then $f_3$ is $4$-sided, and the $D$-lines $d_3$ and $d_4$, and thus $d_2$, are determined.
But, $f_4$ can be neither $3$-sided nor $4$-sided, because
the existence of $b_2$ disturbs an $A$-line and a $B$-line as above.
We have thus shown that $f_1$ is $4$-sided.

Next, we  claim that $f_2$ is $4$-sided.
Assume not.
Then vertex $D$ is located, 
and thus $d_4$ is determined.
Also, $f_3$ is $4$-sided.
If $f_3$ is incident with $x_2$, then
the $D$-line $d_2$ is determined, and thus
$d_3$ cannot be drawn.
Hence $f_3$ is incident with another triple crossing point
as in the first of Figure \ref{fig:caseii-2-2-1a}, where
an $A$- or $C$-line goes through.
If it is a $C$-line, then
$f_4$ cannot be $3$-sided.
Hence it is an $A$-line, in particular, $a_2$.
See the second of Figure \ref{fig:caseii-2-2-1a}.
Then $d_2$ cannot be drawn.
Thus we have specified two $4$-sided faces $f_1$ and $f_2$.

\begin{figure}[tb]
\includegraphics*[scale=0.7]{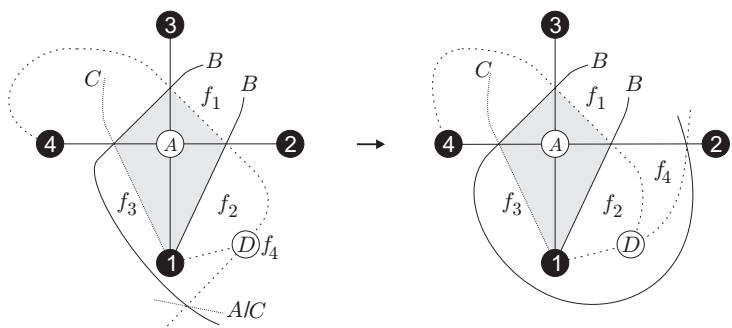}
\caption{The case where $f_2$ is $3$-sided}
\label{fig:caseii-2-2-1a}
\end{figure}

Now, $f_3$ is $3$-sided.
By property ($\ast$), the right lower line of $f_2$ is
an $A$- or $C$-line.
See the first of Figure \ref{fig:caseii-2-2-1a-final0}.
If it is an $A$-line, then it is $a_2$.
But this is impossible, because the right upper line of $f_2$ already meets
$a_2$.
Thus the right lower line of $f_2$ is a $C$-line.
Furthermore, 
if $f_2$ is incident with $x_4$, then
the situation is drawn as in the second of Figure \ref{fig:caseii-2-2-1a-final0}.
Then $b_2$ is determined, and thus vertex $B$ is located.
However, $b_1$ cannot be drawn.

\begin{figure}[tb]
\includegraphics*[scale=0.7]{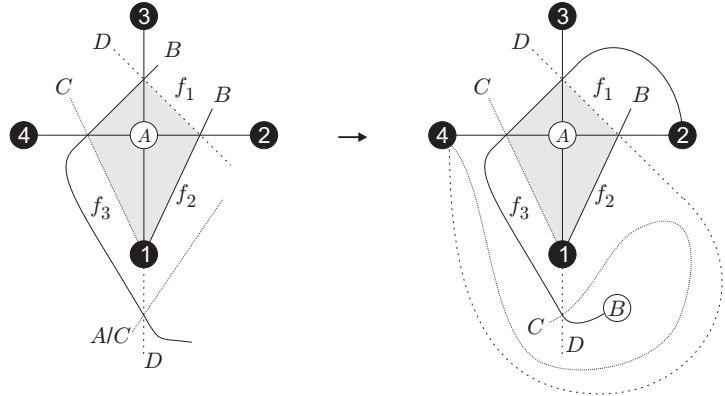}
\caption{$f_1$ and $f_2$ are $4$-sided}
\label{fig:caseii-2-2-1a-final0}
\end{figure}

Thus $f_2$ is incident with a triple crossing point at its right.
See the first of Figure \ref{fig:caseii-2-2-1a-final}.
\begin{figure}[tb]
\includegraphics*[scale=0.7]{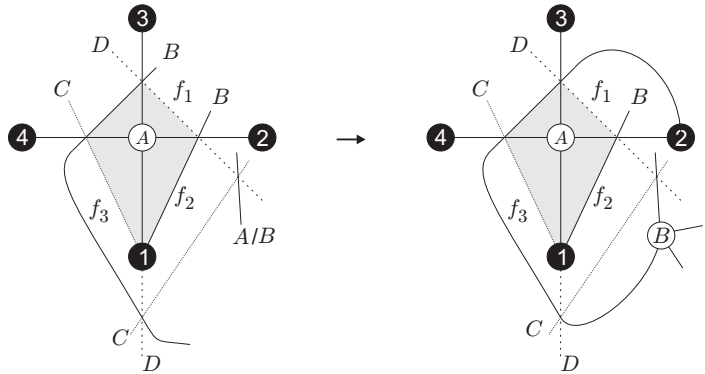}
\caption{$f_1$ and $f_2$ are $4$-sided (continued)}
\label{fig:caseii-2-2-1a-final}
\end{figure}
By property ($\ast$), an $A$- or $B$-line goes through the triple crossing point.
But it cannot be an $A$-line by examining the ($3$-sided) face right above $f_2$.
After locating vertex $B$, the $B$-line $b_2$ is determined.
Then $b_1$ cannot be drawn.
\end{proof}

\begin{claim}
Figure \ref{fig:casei-2}(2) is impossible.
\end{claim}

\begin{proof}
We claim that $f_1$ and $f_4$ are $4$-sided.
Assume that $f_1$ is $3$-sided.
Then vertex $B$ is located, and then the $B$-line $b_2$ is determined.
See the first of Figure \ref{fig:caseii-2-3}.

\begin{figure}[tb]
\includegraphics*[scale=0.7]{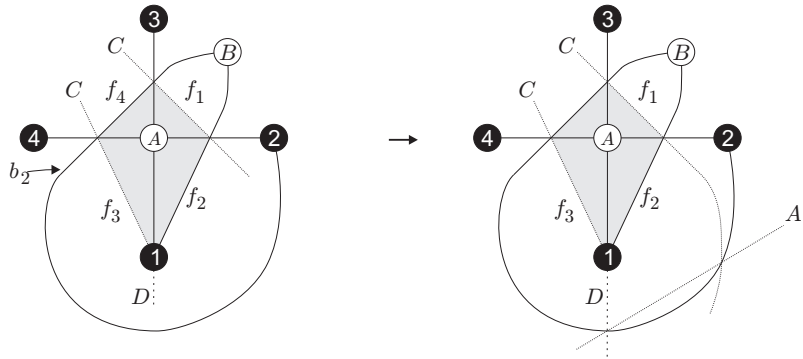}
\caption{The case where $f_1$ is $3$-sided}
\label{fig:caseii-2-3}
\end{figure}

Suppose further that $f_2$ is $4$-sided.
Then $f_2$ is incident with either vertex $C$, or vertex $D$, or
two more triple crossing points. 
If $f_2$ is incident with vertex $C$, then
$c_4$ is determined, and then $c_1$ cannot be drawn.
If $f_2$ is incident with vertex $D$, then
an $A$- or $B$-line appears at the triple crossing point where a $C$-line meets a $D$-line.
But this is impossible by the existence of $b_2$ as in 
the proof of Claim \ref{cl:ii-2-2}.
If $f_2$ is incident with two more triple crossing points, then
a similar argument to the proof
of Claim \ref{cl:ii-2-2} gives a contradiction again (see Figure \ref{fig:caseii-2-3}).
Thus $f_2$ is $3$-sided.

Then the situation is as in the proof of Claim \ref{cl:i-2}, leading to a contradiction.
Thus $f_1$ is $4$-sided.
By the same argument, $f_4$ is $4$-sided.
Then $f_2$ and $f_3$ are $3$-sided.
But this is impossible, because
there are a $B$-line and $C$-line meeting twice.
\end{proof}

This completes the proof of Lemma \ref{lem:ii-2-2}.
\end{proof}

\begin{lemma}\label{lem:ii-2-0}
The number of type II triangles at $A$ is not zero.
\end{lemma}

\begin{proof}
Suppose that there is no type II triangle at $A$.
Then 
the local configuration at $A$ is Figure \ref{fig:casei-3}(1) or (2).

In Figure \ref{fig:casei-3}(1),
if $f_1$ or $f_2$ is $3$-sided, then
we have a contradiction as in the proof of Lemma \ref{lem:i0}.
Thus both are $4$-sided.
By the same reason, $f_3$ and $f_4$ are $4$-sided, a contradiction.

In Figure \ref{fig:casei-3}(2),
if $f_1$ or $f_2$ is $3$-sided, then we have a contradiction as above.
Hence $f_1$ and $f_2$ are $4$-sided.
Thus $f_3$ and $f_4$ are $3$-sided.
Then the $C$-line $c$ meets the $D$-line $d$ twice, a contradiction.
\end{proof}

\begin{proposition}\label{prop:caseii}
Case (2) is impossible.
\end{proposition}

\begin{proof}
This follows from Lemmas 
\ref{lem:ii3}, \ref{lem:ii2}, \ref{lem:ii1}, \ref{lem:ii0}, \ref{lem:ii-2-4}, \ref{lem:ii-2-2} and \ref{lem:ii-2-0}.
\end{proof}

\begin{theorem}\label{thm:k44}
$K_{4,4}$ does not admit a semi-regular drawing.
\end{theorem}

\begin{proof}
This follows from Propositions \ref{prop:casei} and \ref{prop:caseii}.
\end{proof}

\section{$K_{n,3}$}\label{sec:kn3}

Let $G=K_{n,3}$ with $n\ge 5$.
In this section, we show that if $n\ne 6$ then
$G$ does not admit a semi-regular drawing.
Hereafter, we assume that $n\ge 5$ and $n\ne 6$.

\subsection{Exceptional faces}

Let $V_1$ and $V_2=\{A, B, C\}$ be the partite sets of $G$.
As before, we refer to a vertex of $V_1$ (resp. $V_2$) as a black (resp. white) vertex.
Any black vertex is incident with an $A$-line, $B$-line and $C$-line.

Suppose that $G$ admits a semi-regular drawing.
Fix such a drawing, denoted by $G$ again.
Property ($\ast$) holds.  That is, 
at each triple crossing point, an $A$-line, a $B$-line and a $C$-line meet.
Let $k$ be the number of triple crossing points.
Add a new vertex to each triple crossing point.
Then we have a plane graph $G'$ with $n+3+k$ vertices and $3n+3k$ edges.
Since $3(n+3+k)-(3n+3k)-6=3$,
the faces of $G'$ are $3$-sided, except at most three faces, by Lemma \ref{lem:faces}.
We refer to a non-triangular face as an \textit{exceptional face}. 
More precisely, Lemma \ref{lem:faces} claims that
either
\begin{enumerate}
\item $G'$ has only one exceptional face, which is $6$-sided; or
\item $G'$ has just two exceptional faces, which are $5$-sided and $4$-sided, respectively; or
\item $G'$ has just three exceptional faces, which are $4$-sided.
\end{enumerate}

As before, a face of $G$ means that of $G'$.
Let $N$ be the number (counted with multiplicities) of white vertices which are incident with exceptional faces.
Then $0\le N\le 6$, because two white vertices are not adjacent in $G'$.
Since a white vertex is not a cut-vertex of $G'$, 
a white vertex cannot appear around one exceptional face twice.

\begin{lemma}\label{lem:same-pair-exc}
Two exceptional faces are not incident with the same pair of white vertices.
\end{lemma}

\begin{figure}[tb]
\includegraphics*[scale=0.7]{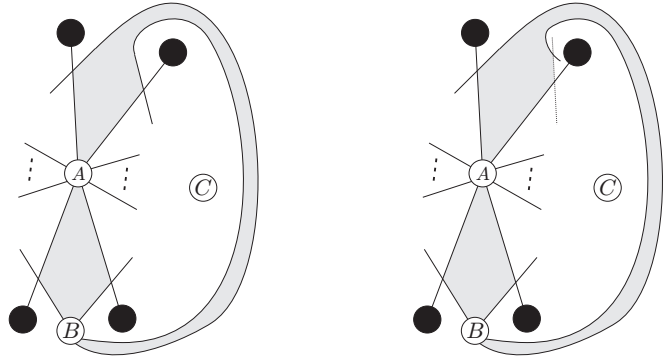}
\caption{Two exceptional faces with the same pair of white vertices}
\label{fig:k3p-N6}
\end{figure}

\begin{proof}
Suppose that two exceptional faces $f$ and $f'$ are incident with white vertices $A$ and $B$, say.
Then both $f$ and $f'$ are $4$-sided, or one is $4$-sided and the other $5$-sided.
See Figure \ref{fig:k3p-N6}.
Recall that any black vertex is incident with a $C$-line.
Thus, in any case, we cannot place $C$-lines. 
\end{proof}

\begin{figure}[tb]
\includegraphics*[scale=0.7]{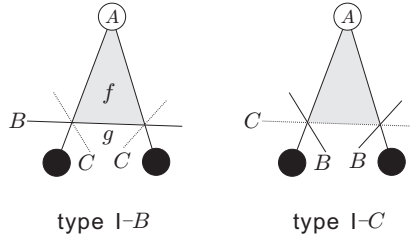}
\caption{A type I-$B$ triangle and a type I-$C$ triangle}
\label{fig:typeI}
\end{figure}

Recall that there are two types of triangles at a white vertex as shown in Figure \ref{fig:triangle}.
Let $X=B$ or $C$.
At vertex $A$, if a type I triangle is bounded by two $A$-lines and an $X$-line,
then it is said to be of \textit{type} I-$X$.  See Figure \ref{fig:typeI}.
Furthermore, a type I-$X$ triangle is said to be \textit{good\/}
if the face sharing the $X$-line with the type I-$X$ triangle is $3$-sided.
Otherwise, it is \textit{bad}.
In particular, a bad type I triangle is adjacent to an exceptional face, which
is referred to as its \textit{associated exceptional face\/}.

\begin{lemma}\label{lem:typeIX-AX}
Let $\{X,Y\}=\{B,C\}$.
If there is a good type I-$X$ triangle at vertex $A$, then
there is neither an exceptional face incident with both $A$ and $Y$, nor
another good type I-$X$ triangle at $A$.
In particular, the number of good type I triangles is at most two.
\end{lemma}

\begin{figure}[tb]
\includegraphics*[scale=0.7]{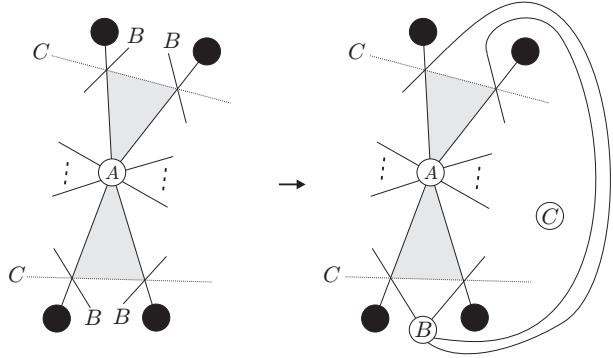}
\caption{Two good type I-$C$ triangles at $A$}
\label{fig:k3p-case1a}
\end{figure}

\begin{proof}
Let $f_1$ be a good type I-$C$ triangle at $A$.
Then the face $f_2$ sharing the $C$-line with $f_1$ is $3$-sided, so
vertex $B$ is located there.
Suppose that an exceptional face $f$ is incident with vertex $A$ and $B$.
Then we have a similar situation to the proof of Lemma \ref{lem:same-pair-exc}.
Hence we cannot place $C$-lines.
The existence of another good type I-$C$ triangle is excluded by a similar argument.
See Figure \ref{fig:k3p-case1a}.
Here, we cannot place the $C$-lines going to the left upper black vertex and the left lower black
vertex simultaneously.
\end{proof}

\begin{lemma}\label{lem:badI}
Suppose that there is a bad type I triangle $f$ at vertex $A$.
Let $g$ be the associated exceptional face of $f$.
If $g$ is $k$-sided ($4\le k\le 6$), then $g$ is incident with at most $k-4$ white vertices.
\end{lemma}

\begin{proof}
We may assume that $f$ is a bad type I-$B$ triangle without loss of generality.
As shown in the first of Figure \ref{fig:typeI}, 
the boundary of $g$ contains a sequence of
\[
\text{a $C$-line, a triple crossing point, a $B$-line, a triple crossing point, a $C$-line.}
\]
The existence of this sequence forces $g$ to admit at most $k-4$ white vertices.
\end{proof}

\begin{lemma}\label{lem:badshare}
Two bad type I triangles at vertex $A$ cannot have the same associated exceptional face.
\end{lemma}

\begin{proof}
Let $f_1$ and $f_2$ be bad type I triangles at $A$ whose associated exceptional faces coincide.
Let $g$ be the common associated exceptional face.
When $g$ is $4$- or $5$-sided,
the situation is as shown in Figure \ref{fig:n4a}, where
labels $B$ and $C$ may be exchanged.
\begin{figure}[tb]
\includegraphics*[scale=0.7]{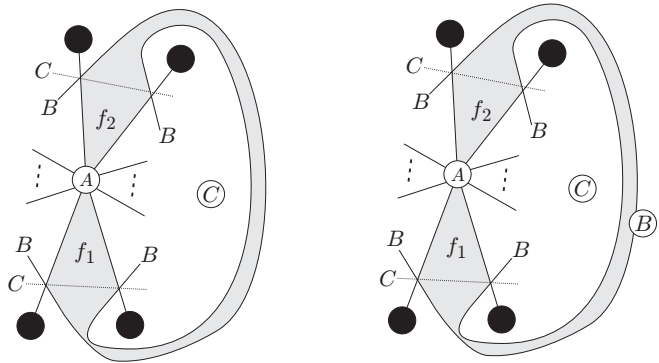}
\caption{Two bad type I triangles with the same associated $4$- or $5$-sided exceptional face}
\label{fig:n4a}
\end{figure}
Then we cannot place $C$-lines as before.
If $g$ is $6$-sided, then the situation is as shown in Figure \ref{fig:n4b}.
Similarly, we cannot draw $C$-lines.
\begin{figure}[tb]
\includegraphics*[scale=0.7]{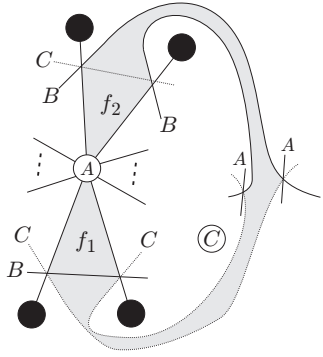}
\caption{Two bad type I triangles with the same associated $6$-sided exceptional face}
\label{fig:n4b}
\end{figure}
\end{proof}

\begin{lemma}\label{lem:badImore}
Suppose that $G$ satisfies one of the following conditions.
\begin{enumerate}
\item $G$ has a $5$-sided exceptional face incident with
exactly one white vertex and a $4$-sided exceptional face incident with at least one white vertex.
\item $G$ has three $4$-sided exceptional faces, only one of which is incident with
no white vertex.
\end{enumerate}
Then there is at most one bad type I triangle at vertex $A$.
\end{lemma}

\begin{proof}
Let $f_i$ be a bad type I triangle at $A$, and let $g_i$ be the associated
exceptional face for $i=1,2$.
By Lemma \ref{lem:badI}, we have $g_1=g_2$.
But this contradicts Lemma \ref{lem:badshare}.
\end{proof}

\begin{lemma}\label{lem:bad3}
There is at most three bad type I triangles at vertex $A$.
\end{lemma}

\begin{proof}
Suppose that there are four bad type I triangles at vertex $A$.
Since $G$ has at most three exceptional faces and
each bad type I triangle is adjacent to an exceptional face,
there exist two bad type I triangles
whose associated exceptional faces coincide, contradicting Lemma \ref{lem:badshare}.
\end{proof}

\begin{figure}[tb]
\includegraphics*[scale=0.7]{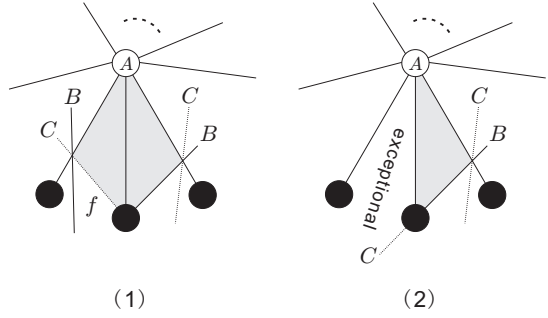}
\caption{Type II triangles at $A$}
\label{fig:k3p-case1b}
\end{figure}

Recall that an adjoint pair of type II triangles at vertex $A$ is
a pair of type II triangles sharing an $A$-line fully.
See Figure \ref{fig:k3p-case1b}(1).

\begin{lemma}\label{lem:typeIIpair}
Suppose that $G$ satisfies one of the following conditions.
\begin{enumerate}
\item $G$ has a single $6$-sided exceptional face, which is incident with all white vertices.
\item $G$ has a $5$-sided exceptional face incident with two white vertices
and a $4$-sided exceptional face incident with a white vertex.
\item $G$ has three $4$-sided exceptional faces, each of which is
incident with a white vertex.
\end{enumerate} 
Then there is no adjoint pair of type II triangles at vertex $A$.
Hence, if there is a type II triangle at $A$, then it shares an $A$-line fully with
an exceptional face (Figure \ref{fig:k3p-case1b}(2)).
\end{lemma}

\begin{proof}
Suppose that there is an adjoint pair of type II triangle at $A$.
Then $f$, indicated in Figure \ref{fig:k3p-case1b}(1), is not $3$-sided.
However, $f$ cannot be an exceptional face from the assumption, a contradiction.
\end{proof}

\begin{lemma}\label{lem:adjointpair}
There is at most one adjoint pairs of type II triangle at vertex $A$.
\end{lemma}

\begin{proof}
As in the proof of Lemma \ref{lem:typeIIpair},
each adjoint pair of type II triangles yields an exceptional face.
If two such pairs share the same exceptional face,
then the exceptional face must be a $6$-sided face
without a white vertex.
The situation is as shown in Figure \ref{fig:k3p-pair}.
\begin{figure}[tb]
\includegraphics*[scale=0.7]{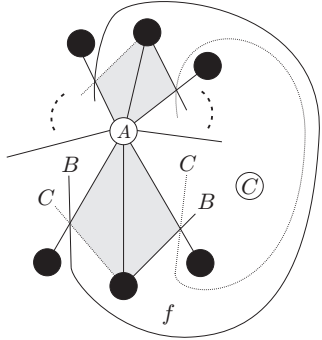}
\caption{Two adjoint pairs of type II triangles sharing the same exceptional face}
\label{fig:k3p-pair}
\end{figure}
Then we cannot draw $C$-lines to the left two black vertices from vertex $C$.

Hence the number of adjoint pairs of type II triangles is no greater than
the number of exceptional faces.
When $G$ has a $6$-sided exceptional face, we have the conclusion.

Assume that $G$ has at least two, then two or three, exceptional faces.
Suppose that there are two adjoint pairs of type II triangle.
Then these pairs correspond to distinct 
exceptional faces as above.
In any case, there exists an adjoint pair of type II triangles which
yields a $4$-sided exceptional face $f$.
Then the situation is as shown in Figure \ref{fig:n2}.
\begin{figure}[tb]
\includegraphics*[scale=0.7]{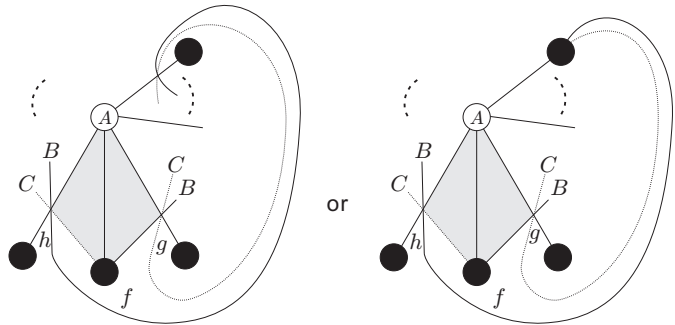}
\caption{An adjoint pair of type II triangles adjacent to a $4$-sided face}
\label{fig:n2}
\end{figure}
This implies that both $g$ and $h$ are exceptional.
Another adjoint pair of type II triangles yields one more exceptional face.
Thus $G$ would have $4$ exceptional faces, a contradiction.
\end{proof}

\begin{lemma}\label{lem:odd}
If $n$ is odd, then each white vertex is incident with an exceptional face, hence $N\ge 3$. 
\end{lemma}

\begin{proof}
Assume that only triangles appear at a white vertex, $A$, say.
Each triangle at $A$ is incident with either a $B$-line or a $C$-line.
Moreover, such triangles appear alternatively around $A$.
Hence $n$ must be even.
\end{proof}

\begin{lemma}\label{lem:even}
Suppose that vertex $A$ is incident with only one exceptional face $f$.
If $f$ is $4$-sided and incident with two white vertices, then $n$ is even.
\end{lemma}

\begin{proof}
We may assume that $f$ is incident with $A$ and $B$.
Then each triangle adjacent to $f$ at $A$ is incident with a $C$-line.
By the same reason as the proof of Lemma \ref{lem:odd},
$n$ is even.
\end{proof}

\subsection{Reduction}

\begin{lemma}\label{lem:N6}
$N\ne 6$.
\end{lemma}

\begin{proof}
Let $N=6$.  This happens only when $G$ has three $4$-sided exceptional faces, each of which
is incident with two white vertices.
In particular, a face incident with two adjacent triple crossing points is $3$-sided.
By Lemma \ref{lem:same-pair-exc}, there is one exceptional face
for each pair of $A$, $B$, $C$.
We examine the local configuration at vertex $A$.

By Lemmas \ref{lem:typeIX-AX} and \ref{lem:badI},
there is no type I triangle.
By Lemma \ref{lem:typeIIpair},
each type II triangle is adjacent to an exceptional face.

If two exceptional faces at $A$ share an $A$-line, then
there are at most two type II triangles.
This implies that $n\le 4$, a contradiction.
(We remark that this situation can happen when $n=4$.)
Otherwise, there are at most four type II triangles.
In fact, both sides of a type II triangle cannot be exceptional faces,
because there is only one exceptional ($4$-sided) face for each pair of $A$, $B$, $C$. 
(See Figure \ref{fig:k3p-case1b}(2).)
Hence there are exactly four type II triangles and two exceptional faces around $A$,
giving $n=6$, a contradiction.
\end{proof}

\begin{lemma}\label{lem:N5}
$N\ne 5$.
\end{lemma}

\begin{proof}
Assume $N=5$.
This happens only when $G$ has three $4$-sided exceptional faces $f_1,f_2,f_3$, two of which
are incident with two white vertices, the other to one white vertex.
By Lemma \ref{lem:same-pair-exc}, we may assume that
$f_1$ is incident with $A$ and $B$, $f_2$ is incident with $B$ and $C$,
and $f_3$ is incident with $B$ or $C$.

We examine the local configuration at $A$.
By Lemma \ref{lem:badI}, there is no bad type I triangle.
Assume that $f_3$ is incident with $B$.
By Lemma \ref{lem:typeIX-AX}, a good type I-$C$ triangle is impossible, and
at most one good type I-$B$ triangle is possible.
By Lemma \ref{lem:typeIIpair},
there are at most two type II triangles, which are adjacent to $f_1$.
Hence we have $n\le 4$, a contradiction.

The case where $f_3$ is incident with $C$ is similar.
\end{proof}

\begin{lemma}\label{lem:N4}
$N\ne 4$.
\end{lemma}

\begin{proof}
Assume $N=4$.
Then $G$ has at least two exceptional faces.

First, suppose that $G$ has a $5$-sided face $f$ and a $4$-sided face $f'$.
Then both $f$ and $f'$ are incident with two white vertices.
We may assume that $f$ are incident with $A$, $B$, and $f'$ to $B$, $C$.
This case is handled by the same argument as in 
the second paragraph of the proof of Lemma \ref{lem:N5}.

Next, suppose that $G$ has three $4$-sided faces $f_1$, $f_2$, $f_3$.
By Lemma \ref{lem:same-pair-exc}, there are five possibilities for 
three exceptional faces as shown in Table \ref{tb:N4}, up to renaming.
By Lemmas \ref{lem:odd} and \ref{lem:even}, $n\ge 8$ in any case.

\begin{table}
\begin{tabular}{|c|c|c|c|}
\hline
 & $f_1$ & $f_2$ & $f_3$ \\
\hline
(a) & $A$, $B$ & $B$, $C$ & none  \\
\hline
(b) & $A$, $B$ & $B$ & $B$  \\
\hline
(c) & $B$, $C$ & $B$ & $C$  \\
\hline
(d) & $A$, $B$ & $C$ & $C$  \\
\hline
(e) & $A$, $B$ & $B$ & $C$  \\
\hline
\end{tabular}
\vskip 10pt
\caption{Five possibilities}\label{tb:N4}
\end{table}

From (b) to (e), 
there is neither bad type I triangle (by Lemma \ref{lem:badI}) nor
adjoint pair of type II triangle (by Lemma \ref{lem:typeIIpair}).
Moreover, there are at most two good type I triangles (by Lemma \ref{lem:typeIX-AX}) and
at most two type II triangles (by Lemma \ref{lem:typeIIpair}).
This gives $n\le 5$, a contradiction.

Consider (a).
By Lemma \ref{lem:typeIX-AX}, there is no good type I-$C$ triangle at $A$, and 
at most one good type I-$B$ triangle is possible.
By Lemma \ref{lem:badImore},
there is at most one bad type I triangle.
In total, there are at most two type I triangles.
By Lemma \ref{lem:adjointpair}, there is at most one adjoint pair of type II triangles, and
further at most two type II triangles, which are adjacent to $f_1$, can be possible.
Hence we have $n\le 7$, a contradiction.
\end{proof}

\begin{lemma}\label{lem:N3}
$N\ne 3$.
\end{lemma}

\begin{proof}
Assume $N=3$.  We divide the proof into three cases, according to the set of exceptional faces of $G'$.

\medskip
\textit{Case 1. $G$ has a single $6$-sided exceptional face.}

Let $f$ be the exceptional face.
Then each white vertex is incident with $f$.
By Lemmas \ref{lem:typeIX-AX} and \ref{lem:badI}, there is no type I triangle.
By Lemma \ref{lem:typeIIpair}, there are at most two type II triangles adjacent to $f$.
Hence $n\le 3$, a contradiction.

\medskip
\textit{Case 2. $G$ has a $5$-sided exceptional face and a $4$-sided exceptional face.}

Let $f_1$ and $f_2$ be the $5$-sided, $4$-sided exceptional faces, respectively.
According to white vertices incident with them, there are four possibilities as in Table \ref{tb:N3}, up to renaming.

\begin{table}
\begin{tabular}{|c||c|c||c|c|c|c|}
\hline
 & $f_1$ & $f_2$ & \text{good I-$B$} & \text{good I-$C$}& \text{bad I} & \text{adjoint II pair}\\
\hline
(a) & $B$, $C$ & $B$ & $\le 1$ & $\le 1$ & $\times$ &$\times$ \\
\hline
(b) & $A$, $B$ & $C$ & $\le 1$ &$\times$ & $\times$ & $\times$ \\
\hline
(c) & $B$ & $B$, $C$ &  $\le 1$ & $\le 1$ & $\le 1$ & $\le 1$ \\
\hline
(d) & $C$ & $A$, $B$ & $\le 1$ & $\times$ & $\le 1$  & $\le 1$ \\
\hline
\end{tabular}  
\vskip 10pt
\caption{Four possibilities and triangles at $A$}\label{tb:N3}
\end{table}

By Lemmas \ref{lem:odd} and \ref{lem:even}, we have $n\ge 8$, except case (b).

(a) There is neither bad type I triangle (by Lemma \ref{lem:badI}) nor type II triangle (by Lemma \ref{lem:typeIIpair}).
By Lemma \ref{lem:typeIX-AX}, there are at most two good type I triangles.
So, $n\le 2$, a contradiction.

(b) By Lemma \ref{lem:typeIIpair}, there are at most two type II triangles, which are incident with $f_1$.
There is neither bad type I triangle nor good type I-$C$ triangle.
Thus $n\le 4$, a contradiction.

(c) There can be a good type I-$X$ triangle at $A$ for $X\in \{B,C\}$.
Hence there are at most two good type I triangles by Lemma \ref{lem:typeIX-AX}.
By Lemma \ref{lem:badImore}, there is at most one bad type I triangle.
There is at most one adjoint pair of type II triangles by Lemma \ref{lem:adjointpair}.

\begin{claim}\label{cl:n3c}
A bad type I triangle does not coexist with an adjoint pair of type II triangles.
\end{claim}

\begin{proof}
Suppose that there is an adjoint pair of type II triangles.
Then
it is adjacent to $f_1$ (see Figure \ref{fig:k3p-case1b}(1)).
If there is a bad type I triangle, then
its associated exceptional face is also $f_1$ by Lemma \ref{lem:badI}.
Hence $f_1$ is not incident with a white vertex, a contradiction.
\end{proof}

In any case, we have $n\le 4$, a contradiction.

(d) There is no good type I-$C$ triangle.
By Lemma \ref{lem:badImore}, there is at most one bad type I triangle.
Claim \ref{cl:n3c} holds again.
There are at most two type II triangles, which are incident with $f_2$, by Lemma \ref{lem:adjointpair}.
In any case, we have $n\le 6$, a contradiction.

\medskip
\textit{Case 3. $G$ has three $4$-sided exceptional faces.}

Let $f_1$, $f_2$, $f_3$ be the exceptional faces.
There are five possibilities as in Table \ref{tb:N3-3}.
Again, we have $n\ge 8$, except case (e), by Lemmas \ref{lem:odd} and \ref{lem:even}.

\begin{table}
\begin{tabular}{|c||c|c|c||c|c|c||c|}
\hline
 & $f_1$ & $f_2$ & $f_3$    &\text{good I} & \text{bad I} & \text{adjoint II pair} & $n$ \\
\hline
(a) & $B$, $C$ & $B$ & none & $\le 2$ & $\le 1$ & $\le 1$ & $\le 5$ \\
\hline
(b) & $A$, $B$ & $C$ & none & $\le 1$ & $\le 1$ & $\le 1$  & $\le 7$ \\
\hline
(c) & $B$ & $B$ &  $B$ & $\le 2$ &$\times$  & $\times$ & $\le 2$ \\
\hline
(d) & $B$ & $B$ & $C$ & $\le 2$  &  $\times$ & $\times$ & $\le 2$ \\
\hline
(e) & $A$ & $B$ & $C$ & $\le 2$& $\times$  & $\times$ & $\le 3$ \\
\hline
\end{tabular}
\vskip 10pt
\caption{Five possibilities and triangles at $A$}\label{tb:N3-3}
\end{table}

The number of good type I triangles is at most two, except (b), by Lemma \ref{lem:typeIX-AX}.
For (b), there is no good type I-$C$ triangle by Lemma \ref{lem:typeIX-AX}.
For (c), (d) and (e), there is neither bad type I triangle (by Lemma \ref{lem:badI}) nor adjoint pair of 
type II triangles (by Lemma \ref{lem:typeIIpair}). 
Thus (c) and (d) are settled as in Table \ref{tb:N3-3}.
For (a) and (b), there is at most one bad type I triangle (by Lemma \ref{lem:badImore}) and
at most one adjoint pair of type II triangles (by Lemma \ref{lem:adjointpair}).
Thus these cases are also settled as in Table \ref{tb:N3-3}.

The remaining case is (e).
If there is no type II triangle, then we have $n\le 3$, a contradiction.
Otherwise, let $g$ be a type II triangle.
Then $g$ is adjacent to the exceptional face $f_1$.
We may assume that $g$ is incident with a $B$-line.
Let $h$ be the face sharing this $B$-line with $g$.
Then $h$ is $3$-sided, so vertex $C$ is located there.
This implies that $f_1$ is incident with $A$ and $C$, a contradiction.
\end{proof}

\begin{lemma}\label{lem:Nge3}
$N\ge 3$.
\end{lemma}

\begin{proof}
Assume $N\le 2$.
We can assume that vertex $A$ is not incident with an exceptional face.
By Lemma \ref{lem:odd}, $n$ is even, so $n\ge 8$.
We estimate the number of triangles at $A$ as before.
By Lemmas \ref{lem:typeIX-AX} and \ref{lem:bad3}, there are at most two good type I triangles
and at most three bad type I triangles.

Since $A$ is not incident with an exceptional face,
type II triangles appear as adjoint pairs.
By Lemma \ref{lem:adjointpair}, there is at most one adjoint pair of type II triangles.
Then we have $n\le 7$, a contradiction.
\end{proof}

\begin{theorem}\label{thm:final}
Let $n\ge 5$ and $n\ne 6$.
Then $K_{n,3}$ cannot admit a semi-regular drawing.
\end{theorem}

\begin{proof}
By Lemma \ref{lem:Nge3}, $N\ge 3$.
However, this is impossible by Lemmas \ref{lem:N6}, \ref{lem:N5}, \ref{lem:N4}, and \ref{lem:N3}.
\end{proof}

\section{Complete bipartite graphs}\label{sec:bipartite}

We have already shown that
$K_{4,4}$, $K_{5,4}$, $K_{5,3}$ and $K_{n,3}$ with
$n\ge 7$ do not admit a semi-regular drawing in Sections \ref{sec:k54}, \ref{sec:k44} and \ref{sec:kn3}.

\begin{theorem}\label{thm:bipartite}
Let $G=K_{n_1,n_2}$.
If $n_2\le 2$, then  $\mathrm{tcr}(G)=0$.
If $n_2\ge 3$, then  $\mathrm{tcr}(G)=\infty$ except
$K_{3,3}$, $K_{4,3}$, $K_{6,3}$, $K_{6,4}$.
Moreover, $\mathrm{tcr}(K_{3,3})=\mathrm{tcr}(K_{4,3})=1$,
$\mathrm{tcr}(K_{6,3})=2$ and $\mathrm{tcr}(K_{6,4})=4$.
\end{theorem}

\begin{proof}
If $n_2\le 2$, then $G$ is planar, and thus $\mathrm{tcr}(G)=0$.
The graph $G$ has $p=n_1+n_2$ vertices and $q=n_1n_2$ edges.
Then $q-3p+6=(n_1-3)(n_2-3)-3$.
Hence if $n_2\ge 5$, or $n_2=4$ and $n_1\ge 7$, then $q-3p+6>0$, and thus $\mathrm{tcr}(G)=\infty$
by Lemma \ref{lem:alg-condition}.

Figure \ref{fig:4partite} (after removing
three edges $v_2v_3$, $v_3v_4$, $v_4v_2$) shows that
$\mathrm{tcr}(K_{3,3})=\mathrm{tcr}(K_{4,3})=1$, since $K_{3,3}$ and
$K_{4,3}$ are not planar.
Note that $\mathrm{tcr}(K_{6,3})\ge 2$, since $\mathrm{cr}(K_{6,3})=6$ (\cite{K})
and $3\mathrm{tcr}(K_{6,3})\ge \mathrm{cr}(K_{6,3})$.
Thus we have that $\mathrm{tcr}(K_{6,3})=2$ from
Figure \ref{fig:6111}
(after removing three edges $v_2v_3$, $v_3v_4$, $v_4v_2$).
Similarly, we have that $\mathrm{tcr}(K_{6,4})=4$ from the fact
$\mathrm{cr}(K_{6,4})=12$ (\cite{K}) and 
Figure \ref{fig:k64}.

\begin{figure}[tb]
\includegraphics*[scale=0.7]{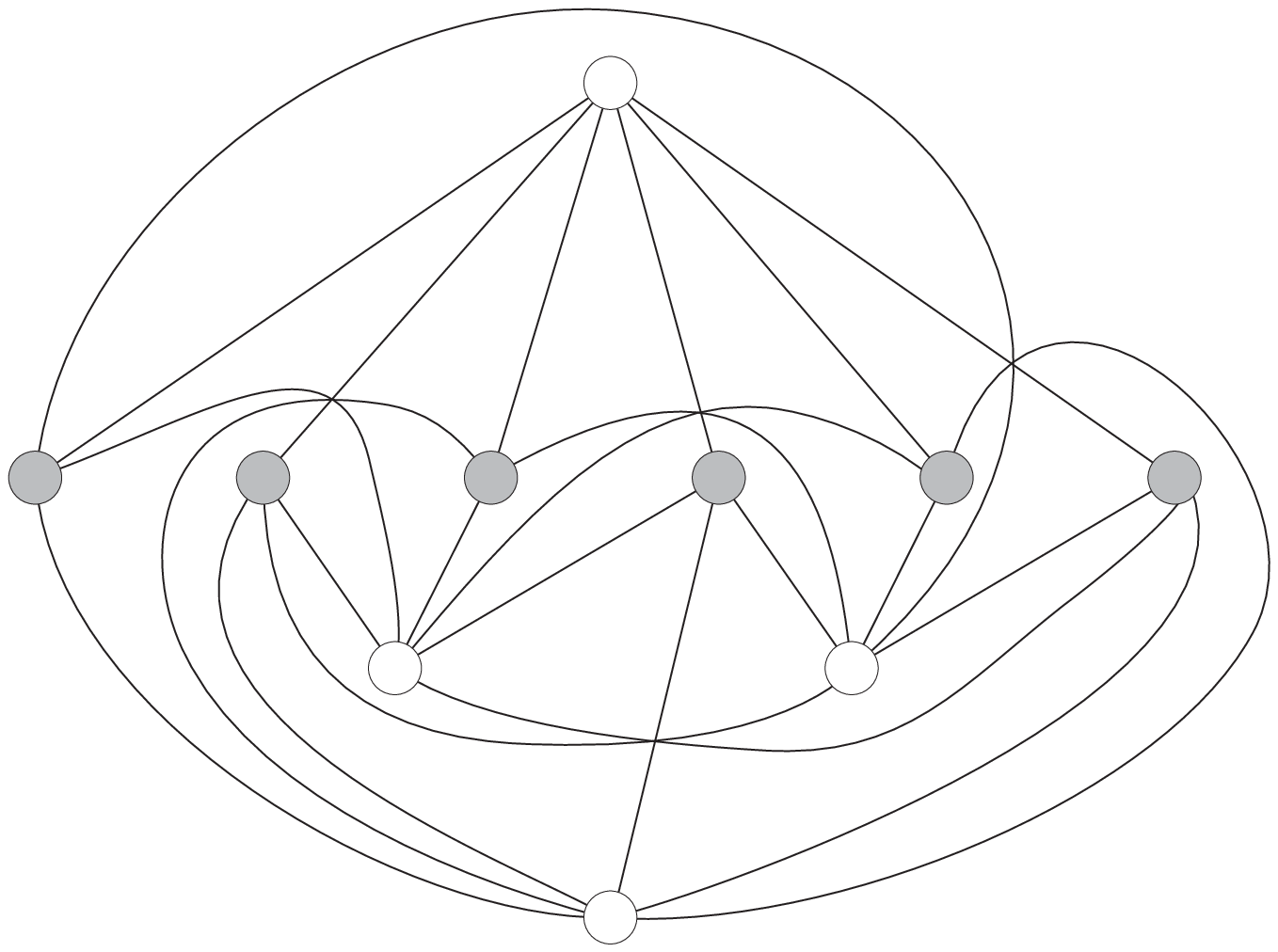}
\caption{$K_{6,4}$}
\label{fig:k64}
\end{figure}

Theorems \ref{thm:k45}, \ref{thm:k44} and \ref{thm:final} show that
$G=K_{n_1,n_2}$ has no semi-regular drawing
for $(n_1,n_2)=(4,4)$, $(5,4)$, $(5,3)$, $(n,3)$ with $n\ge 7$.
\end{proof}

\section{Complete $4$-partite graphs}\label{sec:4partite}

\begin{theorem}\label{thm:4partite}
Let $G=K_{n_1,n_2,n_3,n_4}$.
Then $\mathrm{tcr}(G)=\infty$, except
$K_{n_1,1,1,1}$ with $n_1\in \{1,2,3,4,6\}$.
Also, 
\[
\mathrm{tcr}(K_{n_1,1,1,1})=
\begin{cases}
0 &\text{if $n_1=1,2$},\\
1 &\text{if $n_1=3,4$},\\
2 &\text{if $n_1=6$}.
\end{cases}
\]
\end{theorem}

\begin{proof}
The graph $G$ has $p=\sum_i n_i$ vertices and $q=\sum_{i<j}n_in_j$ edges.
If $n_2\ge 2$, then 
\begin{eqnarray*}
q-3p+6&=& (n_1-1)(n_2-1)+(n_1+n_2)(n_3+n_4-2)+(n_3-3)(n_4-3)-4\\
      &\ge & 1+4(n_3+n_4-2)+(n_3-3)(n_4-3)-4\\
      &=& (n_3+1)(n_4+1)-3\ge 1.
\end{eqnarray*}
Hence we have $\mathrm{tcr}(G)=\infty$ by
Lemma \ref{lem:alg-condition}.

Consider the case where $n_2=n_3=n_4=1$.
If $n_1=1$ or $2$, then $G$ is planar, and thus
$\mathrm{tcr}(G)=0$.
Suppose $n_1\ge 3$.
Then $G$ is non-planar by Lemma \ref{lem:nonplanar}, because $G$ contains
$K_{3,3}$ as a subgraph.
Let $V$ be the partite set of $G$ with $n_1$ elements,
and let $v_2$, $v_3$, $v_4$ be the other vertices of $G$.
Notice that if $G$ admits a semi-regular drawing, then
no edge of the triangle $v_2v_3v_4$ contains
a triple crossing point.
By removing three edges of the triangle from $G$,
we obtain a semi-regular drawing of a complete bipartite
graph $K_{n_1,3}$.
However, this is impossible by Theorem \ref{thm:bipartite},
unless $n_1=3$, $4$ or $6$.
Since 
$K_{3,1,1,1}$ and $K_{4,1,1,1}$ admit a semi-regular drawing with one triple crossing as shown in Figure \ref{fig:4partite},
they have triple crossing number one.
\begin{figure}[tb]
\includegraphics*[scale=0.7]{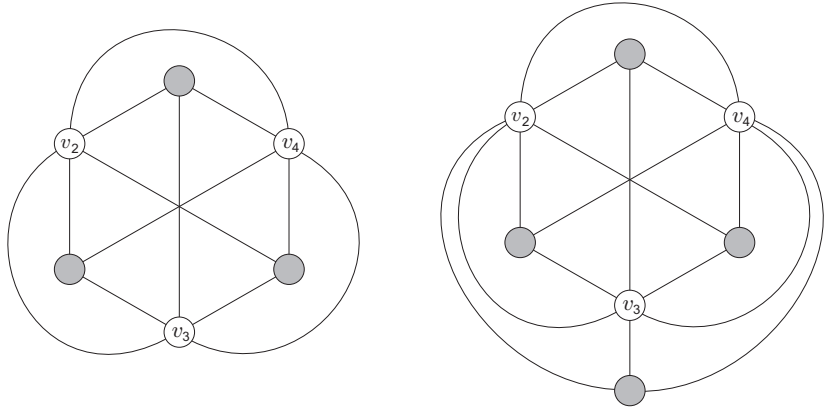}
\caption{$K_{3,1,1,1}$ and $K_{4,1,1,1}$}
\label{fig:4partite}
\end{figure}
Finally, $K_{6,1,1,1}$ admits a semi-regular drawing with two triple crossings
as shown in Figure \ref{fig:6111}.
\begin{figure}[tb]
\includegraphics*[scale=0.7]{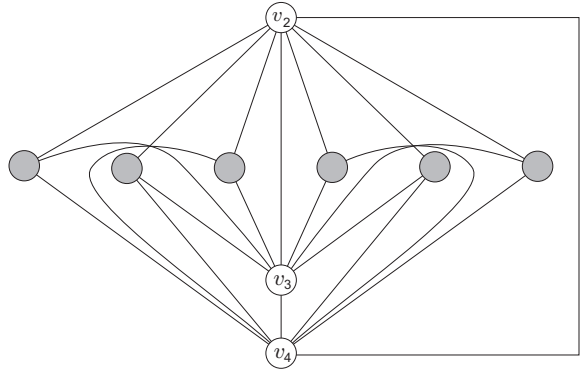}
\caption{$K_{6,1,1,1}$}
\label{fig:6111}
\end{figure}
Since $\mathrm{tcr}(K_{6,3})=2$ by Theorem \ref{thm:bipartite}, $\mathrm{tcr}(K_{6,1,1,1})=2$.
\end{proof}

\section{Complete tripartite graphs}\label{sec:tripartite}

\begin{theorem}\label{thm:tri}
Let $G=K_{n_1,n_2,n_3}$.
Then we have the value of the triple crossing number of $G$
as in Table \ref{table:tri}.
\begin{table}
\begin{tabular}{|c|c|c||c|}
\hline
$n_3$ & $n_2$ & $n_1$ & $\mathrm{tcr}(G)$ \\
\hline\hline
$\ge 3$ & & & $\infty$ \\
\hline
    & $\ge 3$ & & $\infty$ \\
\cline{2-4}
\raisebox{0.1ex}[0cm][0cm]{$2$}
    &  & $\ge 3$ & $\infty$ \\
\cline{3-4}
    &  \raisebox{1ex}[0cm][0cm]{$2$}    &  $2$ & $0$ \\
\hline
\raisebox{-8ex}[0cm][0cm]{$1$}    & $\ge 4$ &  & $\infty$ \\
\cline{2-4}
    &  \raisebox{-1.5ex}[0cm][0cm]{$3$} & $\ge 4$ & $\infty$ \\
\cline{3-4}
    &       & $3$  & $1$ \\
\cline{2-4} 
    &  \raisebox{-5ex}[0cm][0cm]{$2$} & $\ne 2,3,4,6$   & $\infty$ \\
\cline{3-4}
    &       &  $6$  & $2$ \\
\cline{3-4}
    &        &  $3,4$ &  $1$ \\
\cline{3-4}
    &        &  $2$  & $0$ \\
\cline{2-4}
    &  $1$   &   &   $0$ \\
\hline
\end{tabular}
\vskip 10pt
\caption{$\mathrm{tcr}(K_{n_1,n_2,n_3})$}
\label{table:tri}
\end{table}
\end{theorem}

\begin{proof}
The graph $G$ has $p=\sum_i n_i$ vertices and $q=\sum_{i<j}n_in_j$ edges.
Then
\begin{eqnarray*}
q-3p+6 &=& (n_1+n_3-3)(n_2+n_3-3)-n_3^2+3n_3-3\\
       &\ge & (2n_3-3)^2-n_3^2+3n_3-3 \\
       &=& 3(n_3-1)(n_3-2).
\end{eqnarray*}
If $n_3\ge 3$, then $q-3p+6>0$.
Thus we obtain $\mathrm{tcr}(G)=\infty$
by Lemma \ref{lem:alg-condition}.

Let $n_3=2$.
Then $q-3p+6=(n_1-1)(n_2-1)-1$.
If $n_2\ge 3$, then $q-3p+6>0$.
If $n_2=2$, then $q-3p+6>0$, except when $n_1=2$.
For these cases, $\mathrm{tcr}(G)=\infty$ by Lemma \ref{lem:alg-condition} again.
Since $K_{2,2,2}$ is planar, $\mathrm{tcr}(K_{2,2,2})=0$.

Let $n_3=1$.
We have $q-3p+6=(n_1-2)(n_2-2)-1$.
If $n_2\ge 4$, or if $n_2=3$ and $n_1\ge 4$, then $q-3p+6>0$.
For these cases, $\mathrm{tcr}(G)=\infty$.
Since $K_{3,3,1}$ is not planar (it contains $K_{3,3}$ as a subgraph)
and it admits a semi-regular drawing with one triple crossing
as shown in Figure \ref{fig:k331}, we have $\mathrm{tcr}(K_{3,3,1})=1$.

\begin{figure}[tb]
\includegraphics*[scale=0.7]{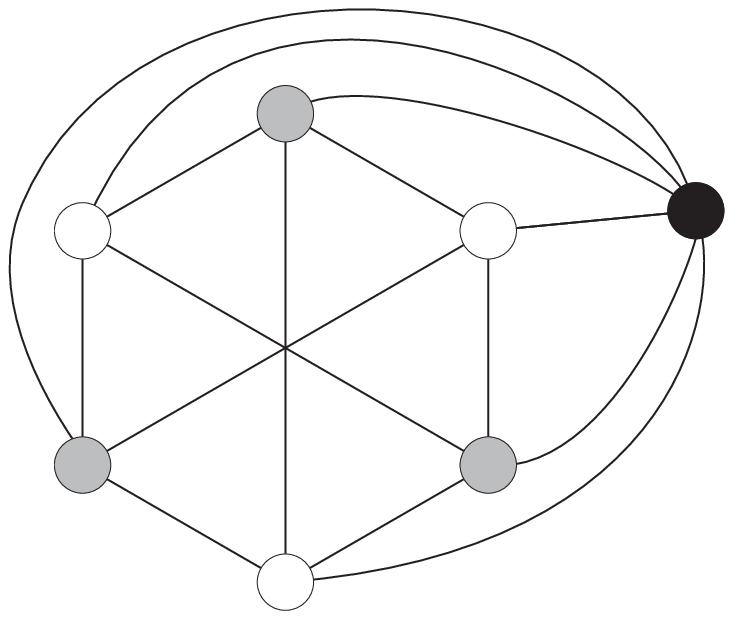}
\caption{$K_{3,3,1}$}
\label{fig:k331}
\end{figure}

The remaining cases are when $n_2=1,2$.
If $n_2=1$ or $n_1=n_2=2$, 
$G$ is planar.
Therefore consider the case where $n_2=2$
and $n_1\ge 3$.
Let $V_1$ and $V_2$ be the partite sets of $G$ with $n_1$ and $n_2$ elements,
respectively, and let $v_3$ be the remaining vertex.
Assume that $G$ admits a semi-regular drawing.
Notice that no edge connecting $v_3$ and a vertex of $V_2$ contains
a triple crossing point, since $|V_2|=2$.
Thus we obtain a semi-regular drawing of a complete bipartite graph $K_{n_1,3}$
by removing two edges between $v_3$ and $V_2$.
Then we have that $n_1=3$, $4$ or $6$ from Theorem \ref{thm:bipartite}.
In either case, $G$ is not planar, since $G$ contains $K_{3,3}$ as a subgraph.
Thus we have that
$\mathrm{tcr}(K_{3,2,1})=\mathrm{tcr}(K_{4,2,1})=1$
by Figure \ref{fig:4partite} (after removing
the edge $v_2v_4$).
Since $\mathrm{tcr}(K_{6,3})=2$ by Theorem \ref{thm:bipartite},
we obtain that $\mathrm{tcr}(K_{6,2,1})=2$ by Figure \ref{fig:6111} (after removing
the edge $v_2v_4$).
\end{proof}


\section{Comments}\label{sec:comment}

In this paper, we require that two edges intersect at most once, and
two edges with a common end-vertex do not intersect.
This is one natural standpoint in the study of the crossing number (\cite{PT0,S}), but
this might be so strong that most complete multipartite graphs
do not admit semi-regular drawings.
If we relax it, then $K_{4,4}$, for example, admits
a semi-regular drawing.

In general, for $n\ge 4$, we can define the $n$-fold crossing number
for a graph $G$ to be the minimal number of $n$-fold crossing points
over all drawings with only $n$-fold crossings.
Clearly, Theorem \ref{thm:5multi} holds for the $n$-fold crossing number.
Furthermore, if $G$ is a non-planar complete $t$-partite graph with $t\ge 3$, then
we can show that $G$ does not admit a drawing with only $n$-fold crossings
by similar arguments to those of Sections \ref{sec:pre}, \ref{sec:4partite} and \ref{sec:tripartite}.
It might be possible to determine the values of this invariant for complete bipartite graphs.

\bibliographystyle{amsplain}

%
%
%
\end{document}